\documentclass[10pt]{amsart}
\usepackage{amssymb}

\theoremstyle{definition}
\newtheorem{thm}{Theorem}[section]
\newtheorem{lem}[thm]{Lemma}
\newtheorem{prp}[thm]{Proposition}
\newtheorem{dfn}[thm]{Definition}
\newtheorem{cor}[thm]{Corollary}
\newtheorem{cnj}[thm]{Conjecture}
\newtheorem{cnv}[thm]{Convention}
\newtheorem{rmk}[thm]{Remark}

\newtheorem{exa}[thm]{Example}

\newtheorem{qst}[thm]{Question}
\newenvironment{pff}{{\emph{Proof:}}}{\QED}

\newcommand{\beq}{\begin{equation}}
\newcommand{\eeq}{\end{equation}}
\newcommand{\beqr}{\begin{eqnarray*}}
\newcommand{\eeqr}{\end{eqnarray*}}
\newcommand{\bal}{\begin{align*}}
\newcommand{\eal}{\end{align*}}
\newcommand{\bei}{\begin{itemize}}
\newcommand{\eei}{\end{itemize}}
\newcommand{\limi}[1]{\lim_{{#1} \to \infty}}

\newcommand{\af}{\alpha}
\newcommand{\bt}{\beta}
\newcommand{\gm}{\gamma}
\newcommand{\dt}{\delta}
\newcommand{\ep}{\varepsilon}
\newcommand{\zt}{\zeta}
\newcommand{\et}{\eta}
\newcommand{\ch}{\chi}

\newcommand{\ld}{\lambda}
\newcommand{\sm}{\sigma}

\newcommand{\ph}{\varphi}

\newcommand{\rh}{\rho}
\newcommand{\om}{\omega}
\newcommand{\ta}{\tau}

\newcommand{\Gm}{\Gamma}

\newcommand{\Sm}{\Sigma}

\newcommand{\Om}{\Omega}

\newcommand{\Q}{{\mathbf{Q}}}
\newcommand{\Z}{{\mathbf{Z}}}
\newcommand{\R}{{\mathbf{R}}}
\newcommand{\C}{{\mathbf{C}}}
\newcommand{\N}{{\mathbf{N}}}

\newcommand{\ev}{{\mathrm{ev}}}

\newcommand{\dist}{{\mathrm{dist}}}
\newcommand{\sa}{{\mathrm{sa}}}
\newcommand{\spec}{{\mathrm{sp}}}

\newcommand{\diag}{{\mathrm{diag}}}
\newcommand{\supp}{{\mathrm{supp}}}
\newcommand{\rank}{{\mathrm{rank}}}

\newcommand{\card}{{\mathrm{card}}}

\newcommand{\dirlim}{\displaystyle \lim_{\longrightarrow}}
\newcommand{\Mi}{M_{\infty}}

\newcommand{\andeqn}{\,\,\,\,\,\, {\mbox{and}} \,\,\,\,\,\,}
\newcommand{\QED}{\rule{0.4em}{2ex}}
\newcommand{\ts}[1]{{\textstyle{#1}}}
\newcommand{\ds}[1]{{\displaystyle{#1}}}
\newcommand{\ssum}[1]{{\ts{ {\ds{\sum}}_{#1} }}}

\newcommand{\ca}{C*-algebra}
\newcommand{\ct}{continuous}
\newcommand{\pj}{projection}

\newcommand{\nbhd}{neighborhood}

\newcommand{\hm}{homomorphism}
\newcommand{\wolog}{without loss of generality}
\newcommand{\Wolog}{Without loss of generality}

\newcommand{\ifo}{if and only if}

\newcommand{\mops}{mutually orthogonal \pj s}

\newcommand{\hme}{homeomorphism}
\newcommand{\mh}{minimal homeomorphism}
\newcommand{\tgca}{transformation group \ca}
\newcommand{\tggp}{transformation group groupoid}

\newcommand{\cfn}{continuous function}
\newcommand{\hsa}{hereditary subalgebra}

\newcommand{\mvnt}{Murray-von Neumann equivalent}
\newcommand{\mvnc}{Murray-von Neumann equivalence}
\newcommand{\KRd}{Kakutani-Rokhlin decomposition}

\renewcommand{\S}{\subset}
\newcommand{\ov}{\overline}
\newcommand{\SM}{\setminus}
\newcommand{\I}{\infty}

\title[Actions of $\Z^d$ on the Cantor
 set]{Crossed products of the Cantor set by free minimal
actions of $\Z^d$}

\author{N.\  Christopher Phillips}

\date{13 May 2002}  

\address{Department of Mathematics, University  of Oregon,
       Eugene OR 97403-1222, USA.}

\email[]{ncp@darkwing.uoregon.edu}

\thanks{Research partially supported by NSF grant DMS 0070776.}

\begin{document}

\begin{abstract}
Let $d$ be a positive integer, let $X$ be the Cantor set, and let
$\Z^d$ act freely and minimally on $X$.
We prove that the crossed product $C^* (\Z^d, X)$ has stable rank one,
real rank zero, and cancellation of \pj s, and that the order on
$K_0 ( C^* (\Z^d, X) )$ is determined by traces.
We obtain the same conclusion for the \ca s of various kinds of
aperiodic tilings.
\end{abstract}

\maketitle

In~\cite{Pt5}, Putnam considered the \ca\  $A$ associated with a
substitution tiling system satisfying certain additional conditions,
and proved that the order on $K_0 (A)$ is determined by the unique
tracial state $\ta$ on $A$.
That is, if $\et \in K_0 (A)$ satisfies $\ta_* (\et) > 0$, then
there is a \pj\  $p \in \Mi (A) = \bigcup_{n = 1}^{\I} M_n (A)$
such that $\et = [p]$.

In this paper, we strengthen Putnam's theorem, obtaining
Blackadar's Second Fundamental
Comparability Question (\cite{Bl4}, 1.3.1) for $A$,
namely that if $p, \, q \in \Mi (A)$ are \pj s such that
$\ta (p) < \ta (q)$ for every tracial state $\ta$ on $A$,
then $p \precsim q$, that is, that $p$ is \mvnt\  to a sub\pj\  of $q$.
We further prove that the \ca\  $A$ has real rank zero~\cite{BP} and
stable rank one~\cite{Rf}.
We also extend the theorem: the same conclusions hold for the \ca s
of some other kinds of aperiodic tilings, and when $A$ is the \tgca\  %
$C^* (\Z^d, X)$ of an
arbitrary free and minimal action of $\Z^d$ on the Cantor set $X$.

We should also mention the recent proof of the gap labelling conjecture
for the Cantor set
(\cite{BBG}, \cite{BO}, \cite{KmP}), which states that the image
of $K_0 (C^* (\Z^d, X))$ under the map to $\R$ induced by a trace is the
subgroup generated by the values of the corresponding invariant measure
on compact open subsets of $X$.

Our results are loosely related to the Bethe-Sommerfeld Conjecture
for quasicrystals in the tight binding approximation.
The tight binding Hamiltonian for a quasicrystal coming from an
aperiodic tiling is a selfadjoint element of the \ca\  of the tiling.
When this \ca\  has real rank zero,
any selfadjoint element has arbitrarily small perturbations
which have finite spectrum, and moreover selfadjoint elements
with totally disconnected spectrum are generic (form a dense
$G_{\dt}$-set) in the set of all selfadjoint elements.

Our proofs are based on the methods of Section~3 of~\cite{Pt5}.
These methods require the presence of a ``large'' AF subalgebra.
For the substitution tilings of~\cite{Pt5}, a suitable subalgebra
is constructed there.
For \tgca s of free minimal actions of $\Z^d$ on the Cantor set,
we obtain this subalgebra by reinterpreting the main result of
Forrest's paper~\cite{Fr} in terms of groupoids.
We actually prove our main results for the reduced \ca s of what we
call almost AF Cantor groupoids.
These form a class of groupoids to which the methods of
Section~3 of~\cite{Pt5} are applicable.
Forrest in effect shows
that the \tggp\  of a free minimal action of $\Z^d$ on the Cantor set
is almost AF, and Putnam in effect shows in Section~2 of~\cite{Pt5}
that the groupoids of the substitution tiling systems considered
there are almost AF.

There are three reasons for presenting our main results in terms of
almost AF Cantor groupoids.
First, the abstraction enables us to focus on just a few key properties.
In particular, as will become clear, for actions of $\Z^d$
we do not need the full strength of the results obtained in~\cite{Fr}.
Second, it seems plausible that groupoids arising in other contexts
might turn out to be almost AF, so that our work would apply elsewhere.
Third, we believe that the methods will work for actions of much more
general discrete groups,
and for actions that are merely essentially free.
Proving this is, we hope, primarily a matter of generalizing
Forrest's construction of \KRd s~\cite{Fr}, and we want to separate the
the details of the generalization from the methods used to obtain
from it results for the crossed product \ca s.

This paper is organized as follows.
The first section presents background material on principal r-discrete
groupoids and their \ca s, especially in the case that the unit space
is the Cantor set, in a form convenient for later use.
In the second section, we define almost AF Cantor groupoids and
present some basic results.
The key technical result, Lemma~\ref{CutToAF}, appears here.
In Sections \ref{Sec:TOKTh}--\ref{Sec:tsr},
we prove results for the reduced \ca\  of
an almost AF Cantor groupoid when the \ca\  is simple,
proving (most of) Blackadar's Second Fundamental Comparability Question
in Section~\ref{Sec:TOKTh}, real rank zero in Section~\ref{Sec:RRZ},
and stable rank one in Section~\ref{Sec:tsr}.
The full statement of
Blackadar's Second Fundamental Comparability Question is obtained
by combining the result of Section~\ref{Sec:TOKTh} with stable rank one.
In Section~\ref{Sec:KRd} we use~\cite{Fr} to show that
free minimal actions of $\Z^d$ on the Cantor set yield
almost AF Cantor groupoids.
The structural results above therefore hold for their \ca s.
In Section~\ref{Sec:QC} we do the same for the groupoids associated
with several kinds of aperiodic tilings, and discuss the relation
to the Bethe-Sommerfeld Conjecture.
The last section contains some open problems and an example
related to the nonsimple case.

I am grateful to the following people for helpful conversations
and email correspondence:
Claire Anantharaman-Delaroche, Jean Bellissard, Tomasz Downarowicz,
Alan Forrest, Johannes Kellendonk,
David Pask, Ian Putnam, Jean Renault, and Christian Skau.
Some of these conversations occurred at the conference on
Operator Algebras and Mathematical Physics at Constan\c{t}a in July
2001, and I am grateful to the organizers of that conference for
support.
Part of this work was done during a visit to the Institute of
Mathematics of the Romanian Academy, and I am grateful to that
institution for its hospitality.

\section{Cantor groupoids}\label{Sec:CG}

\indent
In this section, we fix notation, recall some important definitions,
and establish a few elementary facts.
For groupoid notation and terminology,
we will generally follow Renault's book~\cite{Rn},
with two exceptions.
If $G$ is a groupoid with unit space $G^{(0)}$,
we will refer to its range and source maps $r, \, s : G \to G^{(0)}$,
given by
\[
r (g) = g g^{-1} \andeqn s (g) = g^{-1} g
\]
for $g \in G$.
Also, transformation groups will normally act on the left;
see Example~\ref{CTGEx} below for notation.

We will recall many of the relevant definitions from~\cite{Rn},
since they are scattered through the book.

It is convenient to have a term to describe the basic assumptions
we will be imposing on our groupoids.

\begin{dfn}\label{CGDfn}
A topological groupoid $G$ equipped with a Haar system
(which will be suppressed in the notation) is called a
{\emph{Cantor groupoid}} if the following conditions are satisfied:
\begin{itemize}
\item[(1)]
$G$ is Hausdorff, locally compact, and second countable.
\item[(2)]
The unit space $G^{(0)}$ is compact, totally disconnected, and has
no isolated points (so is homeomorphic to the Cantor set).
\item[(3)]
$G$ is r-discrete in the sense of
Definition~2.6 in Chapter~1 of~\cite{Rn},
that is, $G^{(0)}$ is open in $G$.
\item[(4)]
The Haar system consists of counting measures.
\end{itemize}
\end{dfn}

Using Lemma~\ref{CGElem} below, we can rephrase this in terminology that
has recently become common (see for example~\cite{Hf} and~\cite{Md}) by
saying that a Cantor groupoid is a
second countable locally compact Hausdorff etale
groupoid whose unit space is the Cantor set, equipped with
the Haar system of counting measures.

The following lemma describes some of the immediate properties
of Cantor groupoids.

\begin{lem}\label{CGElem}
Let $G$ be a Cantor groupoid.
Then:
\begin{itemize}
\item[(1)]
For any $x \in G^{(0)}$, the sets $r^{-1} (x), \, s^{-1} (x) \S G$
are discrete.
\item[(2)]
The range and source maps $r, \, s : G \to G^{(0)}$ are local \hme s.
\item[(3)]
$G$ is totally disconnected.
\end{itemize}
\end{lem}

\begin{pff}
(1) This is Lemma~2.7(i) in Chapter~1 of~\cite{Rn}.

(2) This is Lemma~2.7(iii) in Chapter~1 of~\cite{Rn}.

(3)
It follows from Part~(2) that every point $g \in G$ has an open
\nbhd\  which is homeomorphic to an open subset of $G^{(0)}$.
Now use the fact that $G^{(0)}$ is assumed to be totally disconnected.
\end{pff}

\smallskip

We now give the motivating example.

\begin{exa}\label{CTGEx}
Let $X$ be the Cantor set, and let $\Gm$ be a countable discrete
group which acts on $X$.
Then the \tggp\  $\Gm \times X$, equipped with
the Haar system consisting of counting measures, is a Cantor groupoid.

For reference, and to establish conventions,
here are the groupoid operations.
The pairs $(\gm_1, x_1)$ and $(\gm_2, x_2)$ are composable exactly
when $x_1 = \gm_2 x_2$, and then
$(\gm_1, x_1) (\gm_2, x_2) = (\gm_1 \gm_2, \, x_2)$.
The range, source, and inverse are given by
\[
r (\gm, x) = (1, \gm x), \,\,\,\,\,\,
s (\gm, x) = (1, x), \andeqn
(\gm, x)^{-1} = ( \gm^{-1}, \, \gm x).
\]
\end{exa}

\begin{exa}\label{OpenSbGpd}
Let $G$ be a Cantor groupoid, and let $H$ be an open subgroupoid
of $G$ which contains the entire unit space $G^{(0)}$ of $G$
(or, more generally, whose unit space is a compact open subset of
$G^{(0)}$ with no isolated points).
Equip $H$ with the Haar system of counting measures.
Then $H$ is a Cantor groupoid.

To verify this, the only nonobvious condition from
Definition~\ref{CGDfn} is that the
counting measures form a Haar system, and for this the only
issue is the second condition (continuity)
in Definition~2.2 of Chapter~1 of \cite{Rn}.
For $x \in H^{(0)}$,
let $\mu^x$ and $\nu^x$ be the counting measures on
\[
\{ g \in G \colon r (g) = x \} \andeqn
\{ g \in H \colon r (g) = x \}.
\]
Let $f \in C_{\mathrm{c}} (H)$.
We have to show that $x \mapsto \int_H f \, d \nu^x$ is \ct.
Since $H$ is open in $G$, the function $f$ extends to
a function ${\widetilde{f}} \in C_{\mathrm{c}} (G)$
by setting ${\widetilde{f}} (g) = 0$ for $g \in G \setminus H$.
Further,
\[
\int_H f \, d \nu^x = \int_{G} {\widetilde{f}} \, d \mu^x,
\]
which is known to be \ct\  in $x$.
So $H$ is a Cantor groupoid.

As a specific example, we mention the subgroupoids of the
\tggp\  of a \mh\  of the Cantor set implicit in~\cite{Pt1}.
For the groupoid interpretation, see Example~2.6 of~\cite{Pt4}.
\end{exa}

For convenience, we will reformulate several standard definitions
and constructions in our restricted context.

\begin{rmk}\label{CGRmk}
Let $G$ be a Cantor groupoid, or, more generally, a locally
compact r-discrete groupoid with counting measures as the Haar system.

(1)
(See the beginning of Section~1 in Chapter~2 of~\cite{Rn}.)
The product and adjoint in the convolution algebra $C_{\mathrm{c}} (G)$
(the space of \cfn s on $G$ with compact support) are given by:
\[
(f_1 f_2) (g)
   = \ssum{h \in G \colon  r (h) = s (g)} f_1 (g h) f_2 (h^{-1})
\andeqn  f^* (g) = \ov{f (g^{-1})}.
\]
(Note that we write $f_1 f_2$ rather than $f_1 * f_2$.)
The \ca\  $C^* (G)$ is the completion of $C_{\mathrm{c}} (G)$ in a
suitable C* norm; see Definition~1.12 in Chapter~2 of~\cite{Rn}.

(2)
(See Definitions~3.2, 3.4, and~3.12 in Chapter~1 of~\cite{Rn}.)
A Borel measure $\mu$ on $G^{(0)}$ is invariant \ifo\  for every
$f \in C_{\mathrm{c}} (G)$, the numbers
\[
\int_{G^{(0)}}
   \left( \ssum{g \in G \colon  r (g) = x} f (g) \right) \, d \mu (x)
\andeqn
\int_{G^{(0)}}
   \left( \ssum{g \in G \colon  s (g) = x} f (g) \right) \, d \mu (x)
\]
are equal.
(The difference in the expressions is that one sum is over
$r (g) = x$ and the other is over $s (g) = x$.)

(3)
(See the discussion preceding Definition~2.8 in Chapter~2 of~\cite{Rn},
but note that Renault seems to reserve the term
``regular representation'' for the case that the measure $\mu$ is
quasiinvariant.)
Let $\mu$ be a Borel measure on $G^{(0)}$.
Then the regular representation $\pi$ of $C^* (G)$ associated with $\mu$
is constructed as follows.
Define a measure $\nu$ on $G$ by
\[
\int_G f \, d \nu
  = \int_{G^{(0)}}
      \left( \ssum{g \in G \colon  s (g) = x} f (g) \right) \, d \mu (x)
\]
for $f \in C_{\mathrm{c}} (G)$.
(This measure is called $\nu^{-1}$ in~\cite{Rn}.)
Then $\pi$ is the representation on $L^2 (G, \nu)$ determined by the
formula
\[
\langle \pi (f) \xi, \, \et \rangle
  = \int_{G^{(0)}} \left( \ssum{g, h \in G \colon  s (g) = r (h) = x}
      f (g h) \xi (h^{-1}) \ov{\et (g)} \right) \, d \mu (x).
\]

(4)
By comparing formulas, we see that the regular representation as
defined here is the same as
the representation called ${\mathrm{Ind}}_{\mu}$ in~\cite{KS}
before Corollary~2.4.
(See Section~1A of~\cite{KS}.)
\end{rmk}

Because of the way the relevant material is presented in the literature,
and to call attention to the neat formulation of~\cite{KS},
we take some care with the definition of the reduced \ca.

\begin{dfn}\label{RedAlg}
(Section~1A of~\cite{KS}.)
Let $G$ be a locally
compact r-discrete groupoid with counting measures as the Haar system.
We define the reduced \ca\  $C_{\mathrm{r}}^* (G)$ to be the completion
of $C_{\mathrm{c}} (G)$ in the supremum of the seminorms coming from
the representations $\ld_x$, for $x \in G^{(0)}$, defined as follows:
let $G_x = \{ g \in G \colon s (g) = x \}$, and let $C_{\mathrm{c}} (G)$
act on the Hilbert space $l^2 (G_x)$ by
\[
\ld_x (f) \xi (g)
    = \ssum{h \in G \colon  s (h) = x} f (g h^{-1} ) \xi (h).
\]
\end{dfn}

As shown in Theorem~2.3 of~\cite{KS}, this norm comes from a single
canonical regular representation of $C_{\mathrm{c}} (G)$ on a
Hilbert module over $C_0 \left( G^{(0)} \right)$.
This might appropriately be used as the definition of
$C_{\mathrm{r}}^* (G)$.

\begin{lem}\label{DfnAgree}
The reduced \ca\  $C_{\mathrm{r}}^* (G)$ as defined above is the same as
the reduced \ca\  as in Definition~2.8 in Chapter~2 of~\cite{Rn}.
\end{lem}

\begin{pff}
The representation $\ld_x$ used in Definition~\ref{RedAlg} is easily
checked to be the representation of Remark~\ref{CGRmk}(3) coming
from the point mass at $x$.
Given this, and Remark~\ref{CGRmk}(4), the result follows from
Corollary~2.4(b) of~\cite{KS}.
\end{pff}

\smallskip

The following two results are well known (in fact, in greater
generality), but we have been unable to locate references.
The version of the first for full groupoid \ca s and full
crossed products is in~\cite{Rn} (after Definition~1.12 in Chapter~2),
but we have been unable to find a statement of the reduced case.

\begin{prp}\label{CrPisGpC}
Let $X$ be a locally compact Hausdorff space, let the
discrete group $\Gm$ act on $X$, and let $G = \Gm \times X$
be the \tggp\  (as in Example~\ref{CTGEx}).
Then the reduced groupoid \ca\  $C^*_{\mathrm{r}} (G)$ is isomorphic
to the reduced crossed product \ca\  $C^*_{\mathrm{r}} (\Gm, X)$.
\end{prp}

\begin{pff}
Recall that $C^*_{\mathrm{r}} (\Gm, X)$ is the completion
of $C_{\mathrm{c}} (\Gm, \, C (X))$ in a suitable norm
(see 7.6.5 and 7.7.4 of~\cite{Pd}), with multiplication and
adjoint on $C_{\mathrm{c}} (\Gm, \, C (X))$ given by twisted versions
of those in the ordinary group algebra (see 7.6.1 of~\cite{Pd}).
Define
$\ph \colon C_{\mathrm{c}} (G) \to C_{\mathrm{c}} (\Gm, \, C (X))$
by $\ph (f) (\gm) (x) = f (\gm, \, \gm^{-1} x)$.
One checks that $\ph$ is a bijective *-homomorphism.
Moreover, when one identifies $G_{(1, x)}$ with $\Gm$ in the obvious
way,
one finds that $\ph$ transforms the representation $\ld_{(1, x)}$ of
Definition~\ref{RedAlg} into the regular representation $\pi_x$
of $C^* (\Gm, \, C (X))$ determined as in 7.7.1 of~\cite{Pd} by
the point evaluation $\ev_x$, regarded as a one dimensional
representation of $C (X)$.
By Definition~\ref{RedAlg},
the representation $\bigoplus_{x \in X} \ld_{(1, x)}$ is injective on
$C^*_{\mathrm{r}} (G)$, and by Theorem~7.7.5 of~\cite{Pd}
the representation $\bigoplus_{x \in X} \pi_x$ is injective on
$C^*_{\mathrm{r}} (\Gm, X)$.
Therefore $\ph$ determines an isomorphism
$C^*_{\mathrm{r}} (G) \to C^*_{\mathrm{r}} (\Gm, X)$.
\end{pff}

\begin{prp}\label{AlgIncl}
Let $G$ be a locally
compact r-discrete groupoid with counting measures as the Haar system.
Let $H$ be an open subgroupoid, and give $H$ also
the Haar system consisting of counting measures.
(See Example~\ref{OpenSbGpd}.)
Then the inclusion $C_{\mathrm{c}} (H) \subset C_{\mathrm{c}} (G)$
defines an injective \hm\  %
$C_{\mathrm{r}}^* (H) \to C_{\mathrm{r}}^* (G)$.
\end{prp}

\begin{pff}
Let $\ld_x^G$ and $\ld_x^H$ denote the representations of
$C_{\mathrm{c}} (G)$ and $C_{\mathrm{c}} (H)$
used in Definition~\ref{RedAlg}.
It suffices to show that, for every $x \in G^{(0)}$, the representation
$\ld_x^G |_{C_{\mathrm{c}} (H)}$ is a direct sum of representations
$\ld_y^H$ for various $y \in H^{(0)}$, and perhaps a copy of the zero
representation.

Define a relation on $G_x$ by $g \sim h$ exactly when $g h^{-1} \in H$.
Let $Y = \{ g \in G_x \colon g \sim g \}$, which is equal to
$\{ g \in G_x \colon r (g) \in H \}$.
Restricted to $Y$, the relation $\sim$ is an equivalence relation.

Let $E$ be an equivalence class.
One easily checks that the subspace $l^2 (E)$ is invariant
for $\ld_x^G |_{C_{\mathrm{c}} (H)}$.
Choose $g_0 \in E$ and let $y = r (g_0)$.
Then the formula $h \mapsto h g_0$ defines a bijection $H_y \to E$.
(To check surjectivity: if $g \sim g_0$, then $g g_0^{-1}$ is in $H$
and has source $y$, and $g = (g g_0^{-1}) g_0$.)
Further, if we define a unitary $u \colon l^2 (E) \to l^2 (H_y)$
by the formula $(u \xi) (h) = \xi (h g_0)$, then we get
$u \left( \ld_x^G (f) |_H \right) u^* = \ld_y^H (f)$
for all $f \in C_{\mathrm{c}} (H)$.

Finally, one easily checks that $l^2 (G_x \setminus Y)$ is an
invariant subspace for $\ld_x^G |_{C_{\mathrm{c}} (H)}$, and that
the restriction of the representation to this subspace is zero.
\end{pff}

\begin{dfn}\label{GSetDfn}
(Definition~1.10 in Chapter~1 of~\cite{Rn}.)
Let $G$ be a groupoid.
A {\emph{$G$-set}} is a subset $S \S G$ such that the restrictions
of both the range and source maps to $S$ are injective.
\end{dfn}

\begin{lem}\label{GSetCover}
Let $G$ be a Cantor groupoid.
Let $K \S G$ be a compact set.
Then $K$ is a finite disjoint union of compact $G$-sets,
which are open if $K$ is open.
\end{lem}

\begin{pff}
Since $r$ and $s$ are local \hme s, for each $g \in K$ there is
a compact open subset $E (g)$ such that the restrictions of both
$r$ and $s$ to $E (g)$ are injective.
Since $K$ is compact, there are $g_1, \dots, g_n \in K$ such that
$E (g_1), \dots, E (g_n)$ cover $K$.
Then set $K_1 = E (g_1) \cap K$ and, inductively,
\[
K_l = (E (g_l) \SM [E (g_1) \cup \cdots \cup E (g_{l - 1})]) \cap K.
\]
The $K_l$ are compact, and open if $K$ is open,
because the $E (g_l)$ are compact and open,
they are $G$-sets because the $E (g_l)$ are $G$-sets,
and clearly $K$ is the disjoint union of $K_1, \dots, K_n$.
\end{pff}

\begin{lem}\label{GSetNbhd}
Let $G$ be a Cantor groupoid, and let $S \S G$ be a compact $G$-set.
Then there exists a compact open $G$-set $T$ which contains $S$.
\end{lem}

\begin{pff}
Write $S = \bigcap_{n = 1}^{\I} V_n$ for compact open subsets $V_n \S G$
with $V_1 \supset V_2 \supset \cdots \supset S$.
We claim that some $V_n$ is a $G$-set.
Suppose not.
Then there are infinitely many $n$ such that $s |_{V_n}$ is not
injective,
or there are infinitely many $n$ such that $r |_{V_n}$ is not
injective.
We assume the first case.
(The proof is the same for the second case.)
Then for all $n$, the restriction $s |_{V_n}$ is not injective.
Choose $g_n, \, h_n \in V_n$ such that $s (g_n) = s (h_n)$
and $g_n \neq h_n$.
By compactness, we may pass to a subsequence and assume that
$g_n \to g$ and $h_n \to h$.
Then $g, h \in S$ and $s (g) = s (h)$.
If $g \neq h$, we have contradicted the assumption that $S$ is a
$G$-set.
If $g = h$, then every \nbhd\  of $g$ contains two distinct elements,
namely $g_n$ and $h_n$ for sufficiently large $n$, whose images under
$s$ are equal; this contradicts the fact
(Lemma~\ref{CGElem}(2)) that $s$ is a local \hme.
Thus in either case we obtain a contradiction,
so some $V_n$ is a $G$-set.
\end{pff}

\begin{lem}\label{GSetMeas}
Let $G$ be a Cantor groupoid.
Let $\mu$ be a Borel measure on $G^{(0)}$,
and let $\nu$ be the measure on $G$ of Remark~\ref{CGRmk}(3).
Let $L \S G$ be a compact $G$-set.
Then:
\begin{itemize}
\item[(1)]
$\nu (L) = \mu (s (L))$.
\item[(2)]
If $\mu$ is $G$-invariant, then $\nu (L) = \mu (r (L))$.
\end{itemize}
\end{lem}

\begin{pff}
Use Lemma~\ref{GSetNbhd} to choose a compact open $G$-set $V$ which
contains $L$.
Then $s |_V : V \to s (V)$ and $r |_V : V \to r (V)$ are \hme s.
It suffices to prove that if $f : G \to [0, 1]$ is any \cfn\  with
$\supp (f) \S V$ and $f = 1$ on $L$, then
\[
\int_G f \, d \nu = \int_{G^{(0)}} f \circ (s |_V)^{-1} \, d \mu
\]
and, when $\mu$ is $G$-invariant,
\[
\int_G f \, d \nu = \int_{G^{(0)}} f \circ (r |_V)^{-1} \, d \mu.
\]
(We take $f \circ (s |_V)^{-1} = 0$ off $s (V)$ and
$f \circ (r |_V)^{-1} = 0$ off $r (V)$.)
Because $V$ is a $G$-set,
the first equation is just the definition of $\nu$.
For the second, we use invariance of $\mu$ to rewrite
\[
\int_G f \, d \nu
  = \int_{G^{(0)}}
      \left( \ssum{g \in G \colon  r (g) = x} f (g) \right) \, d \mu (x)
\]
(changing the condition $s (g) = x$ in the original sum to
the condition $r (g) = x$).
Now the second equation follows in the same way as the first.
\end{pff}

\smallskip

At this point, we recall some further definitions from~\cite{Rn}.

\begin{dfn}\label{PDfn}
(1)
(Definition~1.1 in Chapter~1 of~\cite{Rn}.)
Let $G$ be a groupoid, and let $x \in G^{(0)}$.
The {\emph{isotropy subgroup}} of $x$ is the set
$\{ g \in G \colon r (g) = s (g) = x \}$.
(It is a group with identity element $x$.)

(2)
(Definition~1.1 in Chapter~1 of~\cite{Rn}.)
A groupoid
$G$ is {\emph{principal}} if every isotropy subgroup is
trivial (has only one element).
Equivalently, whenever $g_1, \, g_2 \in G$
satisfy $r (g_1) = r (g_2)$ and $s (g_1) = s (g_2)$, then $g_1 = g_2$.

(3)
(See Page~35 of~\cite{Rn}.)
Let $G$ be a groupoid.
A subset $E \S G^{(0)}$ is {\emph{invariant}} if whenever $g \in G$
with $s (g) \in E$, then also $r (g) \in E$.

(4)
(Definition~4.3 in Chapter~2 of~\cite{Rn}.)
A locally compact groupoid
$G$ is {\emph{essentially principal}} if for every closed invariant
subset $E \S G^{(0)}$, the set of $x \in E$ with trivial
isotropy subgroup is dense in $E$.
\end{dfn}

The groupoids appearing in the following definition will play a
crucial role in what follows.

\begin{dfn}\label{AFGDfn}
A Cantor groupoid $G$ is called
{\emph{approximately finite}} (AF for short),
if it is the increasing union of a sequence of
compact open principal Cantor subgroupoids, each of which contains the
unit space $G^{(0)}$.
\end{dfn}

In Definition~3.7 of~\cite{GPS2}, and with a weaker condition
(designed to allow unit spaces which are only locally compact),
such a groupoid is called an AF equivalence relation.

The next proposition is included primarily to make the connection
with earlier work.
The corollary will be essential, but it is easily proved directly.

\begin{prp}\label{AFGRmk}
An AF Cantor groupoid is an AF groupoid in the sense of
Definition~1.1 in Chapter~3 of~\cite{Rn}.
An AF groupoid
as defined there is an AF Cantor groupoid \ifo\  its unit space
is compact and has no isolated points.
\end{prp}

\begin{pff}
The first statement follows easily from Lemma~3.4 of~\cite{GPS2}.

The ``only if'' part of the second statement is clear.

To prove the rest, let $G$ be an AF groupoid in the sense of~\cite{Rn},
and assume its unit space is compact and has no isolated points.
By definition, we can write $G$ as the increasing union of a sequence
of open subgroupoids $H_n$, each of which has the same unit space
$G^{(0)}$, and each of which is a disjoint union of a sequence of
elementary groupoids (Definition~1.1 in Chapter~3 of~\cite{Rn})
$H_{n, k}$ of types $N_{n, k} \in \{ 1, 2, \dots, \infty \}$.

Since the unit spaces $H_n^{(0)}$ are compact, for each $n$ there are
only finitely many $H_{n, k}$, that is,
$H_n = \coprod_{k = 1}^{t (n)} H_{n, k}$ with $t (k) < \infty$.
If all $N_{n, k}$ are finite, then each $H_n$ is compact, and we are
done.
Otherwise, write $H_{n, k} = X_{n, k} \times F_{n, k}^2$ with
$\card \left( F_{n, k} \right) = N_{n, k}$, where the groupoid
structure is $(x, r, s) (x, s, t) = (x, r, t)$ and other pairs
are not composable.
Further write $F_{n, k}$ as an increasing union
$F_{n, k} = \bigcup_{d = 1}^{\infty} F_{n, k, d}$, with
$\card \left( F_{n, k, d} \right)
   = \min \left( d, \, \card \left( F_{n, k} \right) \right)$.
Set $H_{n, k, d} = X_{n, k} \times F_{n, k, d}^2$,
which is a subgroupoid of $H_{n, k}$, and set
\[
G_{n, d} = \coprod_{k = 1}^{t (n)} H_{n, k, d},
\]
which is a compact open subgroupoid of $G$ with the same unit
space $G^{(0)}$.
Moreover,
\[
G = \bigcup_{n = 1}^{\infty} \bigcup_{d = 1}^{\infty} G_{n, d}.
\]

We construct inductively a sequence $n \mapsto d (n)$ such that
\[
G_{1, \, d (1)} \subset G_{2, \, d (2)} \subset \cdots
\andeqn
\bigcup_{n = 1}^{\infty} G_{n, \, d (n)} = G.
\]
Take $d (1) = 1$.
Given $d (1), \, d (2), \, \dots, \, d (n)$, note that all
$G_{k, r}$, for $1 \leq k \leq n + 1$ and $r \in \N$, are compact and
contained in the increasing union of open sets
$H_{n + 1} = \bigcup_{r = 1}^{\infty} G_{n + 1, \, r}$.
Therefore we can choose $d (n + 1)$ so large that
$G_{n + 1, \, d (n + 1)}$ contains all of
\[
G_{1, \, n + 1}, \, G_{2, \, n + 1}, \, \dots, \, G_{n + 1, \, n + 1},
\, G_{1, \, d (1)}, \, G_{2, \, d (2)}, \,  \dots, \, G_{n, \, d (n)}.
\]
This gives the desired sequence.

Now set $G_n = G_{n, \, d (n)}$.
\end{pff}

\begin{cor}\label{CStarIsAF}
Let $G$ be an AF Cantor groupoid.
Then $C_{\mathrm{r}}^* (G)$ is an AF algebra.
\end{cor}

\begin{pff}
By Proposition~1.15 in Chapter~3 of~\cite{Rn}, the full
\ca\  $C^* (G)$ is AF.
The reduced \ca\  $C_{\mathrm{r}}^* (G)$ is a quotient
(actually, in this case equal to the full \ca).
\end{pff}

\section{Almost AF groupoids}\label{Sec:AAFG}

\indent
In this section we introduce almost AF Cantor groupoids, and prove
some basic properties.
An almost AF Cantor groupoid contains a ``large'' AF Cantor subgroupoid,
and its reduced \ca\  contains a corresponding ``large'' AF subalgebra.
We establish one to one correspondences between the sets of normalized
traces on the two \ca s, and between them and the sets of
invariant Borel probability measures on the unit spaces of the
two groupoids.
In addition, if the reduced \ca\  of the groupoid is simple, so is
the AF subalgebra.
Moreover, an almost AF Cantor groupoid is essentially principal,
and has an invariant measure whose associated regular representation
is injective on the reduced \ca.

Although our main results involve only almost AF Cantor groupoids
whose reduced \ca s are simple, and we don't know how to
generalize them, it seems worthwhile to attempt to give a definition
which is also appropriate for the nonsimple case.
For more details, see the discussion after Definition~\ref{TAFGDfn}.

\begin{dfn}\label{Thin}
Let $G$ be a Cantor groupoid,
and let $K \S G^{(0)}$ be a compact subset.
We say that $K$ is {\emph{thin}} if for every $n$,
there exist compact $G$-sets $S_1, S_2, \dots, S_n \S G$
such that $s (S_k) = K$ and the sets
$r (S_1), \, r ( S_2), \, \dots, \, r (S_n)$ are pairwise disjoint.
\end{dfn}

\begin{dfn}\label{TAFGDfn}
Let $G$ be a Cantor groupoid.
We say that $G$ is {\emph{almost AF}}
if the following conditions hold:
\begin{itemize}
\item[(1)]
There is an open AF subgroupoid $G_0 \S G$ which contains the
unit space of $G$ and such that whenever $K \S G \SM G_0$ is
a compact set, then $s (K) \S G^{(0)}$ is thin in $G_0$ in the sense of
Definition~\ref{Thin}.
\item[(2)]
For every closed invariant subset $E \S G^{(0)}$, and every nonempty
relatively open subset $U \S E$, there is a $G$-invariant Borel
probability measure $\mu$ on $G^{(0)}$ such that $\mu (U) > 0$.
\end{itemize}
\end{dfn}

This definition is an abstraction of the key ideas in the argument of
Section~3 of~\cite{Pt5}.
Note that $G_0$ is not uniquely determined by $G$.

We will see in Proposition~\ref{SimpleAAF} that condition~(2)
is redundant when $C_{\mathrm{r}}^* (G)$ is simple
(or when $C_{\mathrm{r}}^* (G_0)$ is simple).
In the nonsimple case, we would still like the definition to imply
that $C_{\mathrm{r}}^* (G)$ has stable rank one and real rank zero.
We have three motivations for condition~(2).
First, it seems to be exactly
what is needed to guarantee that the groupoid is essentially principal.
Second, it allows products with a totally disconnected compact
metric space,
regarded as a groupoid in which every element is a unit.
(For a \tggp, this corresponds to forming the product with the trivial
action on such a space.)
Third, there is a free nonminimal
action of $\Z$ on the Cantor set
whose \tggp\  $G$ satisfies condition~(1)
but has an open subset in its unit space which is null for all
$G$-invariant Borel probability measures,
and for which $C_{\mathrm{r}}^* (G)$ does not have stable rank one.
See Example~\ref{NonMinEx}.

\begin{lem}\label{InjRepn}
Let $G$ be a second countable locally compact Hausdorff r-discrete
groupoid with Haar system consisting of counting measures.
Suppose that for every nonempty
open subset $U \S G^{(0)}$, there is a $G$-invariant Borel
probability measure $\mu$ on $G^{(0)}$ such that $\mu (U) > 0$.
Then there exists a $G$-invariant Borel probability measure on $G^{(0)}$
such that the regular representation it determines
(Remark~\ref{CGRmk}(3)) is injective on $C_{\mathrm{r}}^* (G)$.
\end{lem}

\begin{pff}
Let $U_1, U_2, \dots$ form a countable base for the topology of
$G^{(0)}$ consisting of nonempty open sets.
Choose a $G$-invariant Borel
probability measure $\mu_n$ on $G^{(0)}$ such that $\mu_n (U_n) > 0$.
Set $\mu = \sum_{n = 1}^{\infty} 2^{-n} \mu_n$, which is
a $G$-invariant Borel probability measure on $G^{(0)}$
such that $\mu (U_n) > 0$ for all $n$.
Then $\supp (\mu)$ is a closed subset of $G^{(0)}$ such that
$\supp (\mu) \cap U_n \neq \varnothing$ for all $n$.
Therefore $\supp (\mu) = G^{(0)}$.
Now apply Corollary~2.4 of~\cite{KS} and Remark~\ref{CGRmk}(4).
\end{pff}

\begin{cor}\label{AAFInjRepn}
Let $G$ be an almost AF Cantor groupoid.
Then there exists a $G$-invariant Borel probability measure on $G^{(0)}$
such that the regular representation it determines,
as in Remark~\ref{CGRmk}(3), is injective on $C_{\mathrm{r}}^* (G)$.
\end{cor}

We now need a lemma on thin sets.

\begin{lem}\label{MZero}
Let $G$ be a Cantor groupoid,
and let $K \S G^{(0)}$ be a compact subset
which is thin in the sense of Definition~\ref{Thin}.
Then:
\begin{itemize}
\item[(1)]
For every $n$, there exist a compact open set $W$ containing $K$
and compact open $G$-sets $W_1, W_2, \dots, W_n \S G$,
such that $s (W_k) = W$ for all $k$, and such that the sets
$r (W_1), \, r ( W_2), \, \dots, \, r (W_n)$ are pairwise disjoint
compact open subsets of $G^{(0)}$.
\item[(2)]
For every $\ep > 0$, there is a compact open subset $V$ of $G^{(0)}$
such that $K \S V$ and $\mu (V) < \ep$ for every invariant Borel
probability measure $\mu$ on $G^{(0)}$.
\item[(3)]
For every invariant Borel probability measure $\mu$ on $G^{(0)}$,
we have $\mu (K) = 0$.
\end{itemize}
\end{lem}

\begin{pff}
(1)
Using Definition~\ref{Thin}(2),
choose compact $G$-sets $S_1, S_2, \dots, S_n \S G$
such that $s (S_k) = K$ and such that the sets
$r (S_1), \, r ( S_2), \, \dots, \, r (S_n)$ are pairwise disjoint.
Choose disjoint compact open sets $U_1, U_2, \dots, U_n \S G^{(0)}$
such that $r (S_k) \S U_k$.
Use Lemma~\ref{GSetNbhd} to choose compact open
$G$-sets $V_1, V_2, \dots, V_n \S G$ such that $S_k \S V_k$.
Replacing $V_k$ by $V_k \cap r^{-1} (U_k)$, we may assume that
$r (V_1), \, r ( V_2), \, \dots, \, r (V_n)$ are pairwise disjoint.
Since $s$ is a local \hme, the sets
$s (V_1), \, s ( V_2), \, \dots, \, s (V_n)$ are all compact open sets
containing $K$.
Define
\[
W = s (V_1) \cap s (V_2) \cap \cdots \cap s (V_n) \andeqn
W_k = V_k \cap s^{-1} \left( W \right).
\]
Then the $W_k$ are compact open $G$-sets such that $S_k \subset W_k$
for all $k$, such that
$r (W_1), \, r ( W_2), \, \dots, \, r (W_n)$ are pairwise disjoint,
and such that $s (W_k) = W$.

(2)
Let $\ep > 0$.
Choose $n \in \N$ with $\frac{1}{n} < \ep$.
Let $W \S G^{(0)}$ and $W_1, W_2, \dots, W_n \S G$ be as in Part~(1).
Let $\mu$ be any invariant Borel probability measure on $G^{(0)}$.
Let $\nu$ be the measure in Remark~\ref{CGRmk}(3).
By Lemma~\ref{GSetMeas},
\[
\mu (r (W_k)) = \nu (W_k) = \mu (s (W_k)) = \mu (W)
\]
for all $k$.
Since the $r (W_k)$ are disjoint and $\mu \left( G^{(0)} \right) = 1$,
it follows that $\mu (W) \leq \frac{1}{n} < \ep$.

(3)
This is immediate from Part~(2).
\end{pff}

\begin{lem}\label{EssP}
Let $G$ be an almost AF Cantor groupoid.
Then $G$ is essentially principal (Definition~\ref{PDfn}(4)).
\end{lem}

\begin{pff}
Let $E \S G^{(0)}$ be a closed $G$-invariant subset.

Let $G_0$ be as in Definition~\ref{TAFGDfn}(1).
Note that $G_0$ is principal.
If $x \in G^{(0)}$ has nontrivial isotropy, then there is
$g \in G$ with $g \neq x$ such that $r (g) = s (g) = x$.
So $g \not\in G_0$, whence $x \in s (G \SM G_0)$.

Now $G \SM G_0$ is a closed subset of a locally compact
second countable Hausdorff space, and therefore is a countable
union of compact subsets: $G \SM G_0 = \bigcup_{n = 1}^{\I} K_n$.
Each $s (K_n)$ is thin relative to $G^{(0)}$.
Let $U_n$ be the interior of $s (K_n) \cap E$ relative to $E$.
Then Lemma~\ref{MZero}(3) implies that $\mu (U_n) = 0$ for
every $G_0$-invariant Borel probability measure $\mu$ on $G^{(0)}$,
and hence for every $G$-invariant Borel probability measure
$\mu$ on $G^{(0)}$.
So $U_n = \varnothing$ by Definition~\ref{TAFGDfn}(2).
Thus $s (K_n) \cap E$ is nowhere dense in $E$.
It follows that
$s (G \SM G_0) \cap E$ is meager in $E$, and in particular that its
complement is dense in $E$.
So the points in $E$ with trivial isotropy are dense in $E$.
\end{pff}

\smallskip

Now we start work toward the correspondences between the sets of
invariant measures.
The following lemma is the key technical result for taking advantage
of the structure of an almost AF Cantor groupoid, not only here but
in later sections as well.
The main part is (3), in which the products are in the \ca\  of the
AF subgroupoid.
The other parts are given for easy reference.

\begin{lem}\label{CutToAF}
Let $G$ be a Cantor groupoid, and let $G_0$ be an AF subgroupoid
satisfying Part~(1) of the definition of an almost AF Cantor groupoid
(Definition~\ref{TAFGDfn}).
Let $f \in C_{\mathrm{c}} (G)$.
Let $K$ and $L$ be compact open subsets of $G^{(0)}$ such that
\[
K \cap s (\supp (f) \cap [G \SM G_0]) = \varnothing
\andeqn
L \cap r (\supp (f) \cap [G \SM G_0]) = \varnothing.
\]
Let $p = \ch_K$ and $q = \ch_L$.
Then (with convolution products evaluated in $C_{\mathrm{c}} (G)$),
we have:
\begin{itemize}
\item[(1)]
$p, \, q \in C_{\mathrm{c}} (G)$ are \pj s.
\item[(2)]
\[
(f p) (g) = \left\{ \begin{array}{ll}
     f (g)     & \hspace{0.3em}  s (g) \in K \\
     0         & \hspace{0.3em}  s (g) \not\in K
    \end{array} \right.
  \,\,\,\, \andeqn  \,\,\,\,
(q f) (g) = \left\{ \begin{array}{ll}
     f (g)     & \hspace{0.3em}  r (g) \in L \\
     0         & \hspace{0.3em}  r (g) \not\in L
    \end{array} \right..
\]
\item[(3)]
$f p, \, q f \in C_{\mathrm{c}} (G_0)$.
\end{itemize}
\end{lem}

\begin{pff}
Part~(1) is obvious.

To prove Parts~(2) and~(3) for $f p$,
we evaluate $(f \ch_K) (g)$ following Remark~\ref{CGRmk}(1).
We have $\ch_K (h) = 0$ for $h \not\in G^{(0)}$, so the formula
reduces to
\[
(f \ch_K) (g) = f (g) \ch_K (s (g))
\]
for $g \in G$.
This is the formula for $f p$ in Part~(2).
Now suppose $g \in G \SM G_0$.
If $f (g) \neq 0$, then $s (g) \in s (\supp (f) \cap [G \SM G_0])$,
so $\ch_K (s (g)) = 0$.
Thus $g \in G \SM G_0$ implies $(f \ch_K) (g) = 0$.
Certainly $\supp ( f \ch_K )$ is compact, so
$f \ch_K \in C_{\mathrm{c}} (G_0)$, which is Part~(3).

The proof of Parts~(2) and~(3) for $q f$ is similar,
or can be obtained from the case already done by applying it to
$f^*$ and taking adjoints.
\end{pff}

\begin{lem}\label{InvM}
Let $G$ be a Cantor groupoid, and let $G_0$ be an AF subgroupoid
satisfying Part~(1) of the definition of an almost AF Cantor groupoid
(Definition~\ref{TAFGDfn}).
Then every $G_0$-invariant Borel probability measure
on $G^{(0)}$ is $G$-invariant.
\end{lem}

\begin{pff}
Let $\mu$ be a $G_0$-invariant probability measure on $G^{(0)}$.
By assumption, we have
\[
\int_{G^{(0)}}
   \left( \ssum{g \in G \colon  r (g) = x} f (g) \right) \, d \mu (x)
= \int_{G^{(0)}}
   \left( \ssum{g \in G \colon  s (g) = x} f (g) \right) \, d \mu (x)
\]
for all $f \in C_{\mathrm{c}} (G_0)$.
We need to verify this equation for all $f \in C_{\mathrm{c}} (G)$.
It suffices to do this for nonnegative functions $f$.

Let $f \in C_{\mathrm{c}} (G)$ be nonnegative, and let $\ep > 0$.
By Lemma~\ref{GSetCover}, we can write $\supp (f)$ as the disjoint
union of finitely many compact $G$-sets, say $N$ of them.
It follows that for any $x \in G^{(0)}$, we have
\[
\card ( \{g \in \supp (f) \colon  r (g) = x \}) \leq N \andeqn
\card ( \{g \in \supp (f) \colon  s (g) = x \}) \leq N.
\]
Set
\[
K_1 = r (\supp (f) \cap [G \SM G_0]) \andeqn
K_2 = s (\supp (f) \cap [G \SM G_0]).
\]
Then $K_1$ and $K_2$ are thin subsets of $G^{(0)}$,
so Lemma~\ref{MZero}(2)
provides compact open subsets $V_1, \, V_2 \S G^{(0)}$ such that
\[
K_j \S V_j \andeqn \mu (V_j) < \frac{\ep}{N \| f \|_{\infty} + 1}
\]
for $j = 1, 2$.
Define $f_1 = \ch_{G^{(0)} \SM V_1} f$.
By Lemma~\ref{CutToAF}(2),
\[
f_1 (g) = \left\{ \begin{array}{ll}
     0       & \hspace{3em}  r (g) \in V_1 \\
     f (g)   & \hspace{3em}  {\mbox{otherwise}}
    \end{array} \right..
\]
Therefore
\begin{align*}
& \left|
\int_{G^{(0)}}
   \left( \ssum{g \in G \colon  r (g) = x} f (g) \right) \, d \mu (x)
 - \int_{G^{(0)}}
   \left( \ssum{g \in G \colon  r (g) = x} f_1 (g) \right) \, d \mu (x)
     \right|   \\
& \hspace{3em}
 = \int_{V_1}
   \left( \ssum{g \in G \colon  r (g) = x} f (g) \right) \, d \mu (x) \\
& \hspace{3em}
 \leq \mu (V_1) \| f \|_{\I}
    \left[ \sup_{x \in G^{(0)}} \card ( \{g \in G \colon  r (g) = x \})
              \right]
 \leq \mu (V_1) \| f \|_{\I} N < \ep.
\end{align*}
In the following calculation, for the first step we use the previous
estimate, for the second we use
$f_1 \in C_{\mathrm{c}} (G_0)$ (which holds by Lemma~\ref{CutToAF}(3))
and the $G_0$-invariance of $\mu$, and for the third we use
$f_1 \leq f$ (which holds because $f$ is nonnegative):
\begin{align*}
\int_{G^{(0)}}
   \left( \ssum{g \in G \colon  r (g) = x} f (g) \right) \, d \mu (x)
  & < \ep + \int_{G^{(0)}}
   \left( \ssum{g \in G \colon  r (g) = x} f_1 (g) \right) \, d \mu (x)
                    \\
  & = \ep + \int_{G^{(0)}}
   \left( \ssum{g \in G \colon  s (g) = x} f_1 (g) \right) \, d \mu (x)
                    \\
  & \leq \ep + \int_{G^{(0)}}
   \left( \ssum{g \in G \colon  s (g) = x} f (g) \right) \, d \mu (x).
\end{align*}
A similar argument, using $f_2 = f \ch_{G^{(0)} \SM V_2}$ in place
of $f_1$, gives the same inequality with the range and source maps
exchanged.
Since $\ep > 0$ is arbitrary, we get
\[
\int_{G^{(0)}}
   \left( \ssum{g \in G \colon  s (g) = x} f (g) \right) \, d \mu (x)
 = \int_{G^{(0)}}
   \left( \ssum{g \in G \colon  r (g) = x} f (g) \right) \, d \mu (x),
\]
as desired.
\end{pff}

\begin{lem}\label{GpdTr}
Let $G$ be a locally
compact r-discrete groupoid with counting measures as the Haar system.

(1) Let $\mu$ be an invariant Borel probability measure on $G^{(0)}$.
Then the formula
\[
\ta (f) = \int_{G^{(0)}} \left( f |_{G^{(0)}} \right) \, d \mu,
\]
for $f \in C_{\mathrm{c}} (G)$, defines a normalized trace on the
\ca\  $C_{\mathrm{r}}^* (G)$.
Moreover, the assignment $\mu \mapsto \ta$ is injective.

(2) Suppose in addition that $G$ is principal.
Then every normalized trace on $C_{\mathrm{r}}^* (G)$
is obtained from an invariant Borel probability measure $\mu$
on $G^{(0)}$ as in~(1).
\end{lem}

\begin{pff}
This is a special case of Proposition~5.4 in Chapter~2 of~\cite{Rn}.
\end{pff}

\begin{lem}\label{TrAgree}
Let $G$ be a Cantor groupoid, and let $G_0$ be an AF subgroupoid
satisfying Part~(1) of the definition of an almost AF Cantor groupoid
(Definition~\ref{TAFGDfn}).
Let $\ta_1$ and $\ta_2$ be two normalized traces on
$C_{\mathrm{r}}^* (G)$ whose restrictions to $C_{\mathrm{r}}^* (G_0)$
are equal.
Then $\ta_1 = \ta_2$.
\end{lem}

\begin{pff}
Let $f \in C_{\mathrm{c}} (G)$, and let $\ep > 0$.
We prove that $| \ta_1 (f^* f) - \ta_2 (f^* f) | < \ep$.
Such elements are dense in the positive elements of
$C_{\mathrm{r}}^* (G)$, so their linear span is dense in
$C_{\mathrm{r}}^* (G)$, and the result will follow.

\Wolog\  $\| f \| \leq 1$.
Let
\[
K = r (\supp (f^* f) \cap [G \SM G_0]).
\]
Since $f^* f \in C_{\mathrm{c}} (G)$,
the set $K$ is thin in $G_0$ (Definition~\ref{Thin}).
Also, since $f^* f$ is selfadjoint, we have
$s (\supp (f^* f) \cap [G \SM G_0]) = K$.
Choose $n \in \N$ with $n > 2 \ep^{-1}$.  
Use Lemma~\ref{MZero}(1)
to choose a compact open set $W$ containing $K$,
and compact open $G$-sets $W_1, W_2, \dots, W_n \S G$
such that $s (W_k) = W$ and the sets
$r (W_1), \, r ( W_2), \, \dots, \, r (W_n)$ are pairwise disjoint
compact open subsets of $G^{(0)}$.
Let $p = \ch_W$, which is a \pj\  in $C_{\mathrm{c}} (G)$.
The function $v_k = \ch_{W_k}$
defines an element of $C_{\mathrm{c}} (G)$
such that $v_k^* v_k = p$ and $v_k v_k^* = \ch_{r (W_k)}$.
Since the \pj s
$\ch_{r (W_1)}, \, \ch_{r ( W_2)}, \, \dots, \, \ch_{r (W_n)}$
are pairwise orthogonal, it follows that
\[
\ta_1 (p), \, \ta_2 (p) \leq \ts{ \frac{1}{n} }
   < \ts{ \frac{1}{2} } \ep.
\]

By Lemma~\ref{CutToAF}(3), the products $(1 - p) f^* f$
and $f^* f (1 - p)$ are in $C_{\mathrm{c}} (G_0)$.
Since $p \in C_{\mathrm{c}} (G_0)$, it follows that
\[
f^* f - p f^* f p
  = (1 - p) f^* f + p f^* f (1 - p) \in C_{\mathrm{c}} (G_0).
\]
Therefore
\[
\ta_1 (f^* f - p f^* f p) = \ta_2 (f^* f - p f^* f p).
\]
On the other hand,
\[
p f^* f p \leq \| f \|^2 p \leq p,
\]
so
\[
0 \leq \ta_1 (p f^* f p) \leq \ta_1 (p)
   < {\textstyle{\frac{1}{2}}} \ep
 \andeqn
0 \leq \ta_2 (p f^* f p) \leq \ta_2 (p)
   < {\textstyle{\frac{1}{2}}} \ep.
\]
It follows that
\[
| \ta_1 (f^* f) - \ta_2 (f^* f) |
  = | \ta_1 (p f^* f p) - \ta_2 (p f^* f p) | < \ep,
\]
as desired.
\end{pff}

\begin{prp}\label{TrAndMeas}
Let $G$ be a Cantor groupoid, and let $G_0$ be an AF subgroupoid
satisfying Part~(1) of the definition of an almost AF Cantor groupoid
(Definition~\ref{TAFGDfn}).
Then the following sets can all be canonically identified:
\begin{itemize}
\item
The space $M$ of $G$-invariant Borel probability measures on $G^{(0)}$.
\item
The space $M_0$ of $G_0$-invariant
Borel probability measures on $G^{(0)}$.
\item
$T (C_{\mathrm{r}}^* (G))$,
the space of normalized traces on $C_{\mathrm{r}}^* (G)$.
\item
$T (C_{\mathrm{r}}^* (G_0))$,
the space of normalized traces on $C_{\mathrm{r}}^* (G_0)$.
\end{itemize}
The map from $M$ to $M_0$ is the identity.
(Both are sets of measures on $G^{(0)}$.)
The map from $T (C_{\mathrm{r}}^* (G))$ to $T (C_{\mathrm{r}}^* (G_0))$
is restriction of traces (using Lemma~\ref{AlgIncl}).
The maps from $M$ to $T (C_{\mathrm{r}}^* (G))$ and from
$M_0$ to $T (C_{\mathrm{r}}^* (G_0))$ are as in Lemma~\ref{GpdTr}.
\end{prp}

\begin{pff}
The map from $M_0$ to $T (C_{\mathrm{r}}^* (G_0))$ is bijective by
Lemma~\ref{GpdTr}, because $G_0$ is principal.
The map from $M$ to $M_0$ is well defined because $G$-invariant measures
are obviously $G_0$-invariant, and is then trivially injective.
It is surjective by Lemma~\ref{InvM}.
The map from $T (C_{\mathrm{r}}^* (G))$ to $T (C_{\mathrm{r}}^* (G_0))$
is injective by Lemma~\ref{TrAgree}, and
the map from $M$ to $T (C_{\mathrm{r}}^* (G))$
is injective by Lemma~\ref{GpdTr}.
The composite $M \to T (C_{\mathrm{r}}^* (G_0))$ is bijective
by what we have already done, so both these maps must be bijective.
\end{pff}

\smallskip

It is apparently not known whether Lemma~\ref{GpdTr}(2) can be
generalized to essentially principal groupoids, but this proposition
shows that its conclusion is valid for almost AF Cantor groupoids.

We next show how to simplify the verification that a groupoid
is almost AF when its reduced \ca\  is simple.
We need the following well known lemma.
We have been unable to find a suitable reference in the literature,
so we include a proof for completeness.

\begin{lem}\label{UAFHasTr}
Every unital AF algebra $B$ has a normalized trace.
\end{lem}

\begin{pff}
Write $B = \dirlim B_n$ for finite dimensional \ca s $B_n$
and injective unital \hm s $B_n \to B_{n + 1}$.
Let $\ta_n$ be a normalized trace on $B_n$,
and use the Hahn-Banach Theorem to extend
$\ta_n$ to a state $\om_n$ on $B$.
Use Alaoglu's Theorem to find a weak* limit point $\ta$ of the
sequence $(\om_n)$.
It is easily checked that $\ta$ is a trace on $B$.
\end{pff}

\begin{prp}\label{SimpleAAF}
Let $G$ be a Cantor groupoid, and let $G_0$ be an AF subgroupoid
satisfying Part~(1) of the definition of an almost AF Cantor groupoid
(Definition~\ref{TAFGDfn}).
Assume that $C_{\mathrm{r}}^* (G)$ is simple, or that
$C_{\mathrm{r}}^* (G_0)$ is simple.
Then $G$ is an almost AF Cantor groupoid.
\end{prp}

\begin{pff}
We must verify Part~(2) of Definition~\ref{TAFGDfn}.
First, if $C_{\mathrm{r}}^* (G)$ is simple,
Proposition~4.5(i) in Chapter~2 of~\cite{Rn} implies that
there are no nontrivial closed $G$-invariant subsets in $G^{(0)}$.
If $C_{\mathrm{r}}^* (G_0)$ is simple, then for the same reason there
are no nontrivial closed $G_0$-invariant subsets in $G^{(0)}$, and
so certainly no nontrivial closed $G$-invariant subsets in $G^{(0)}$.
Therefore, in either case, it suffices to find a $G$-invariant Borel
probability measure $\mu$ on $G^{(0)}$ such that $\mu (U) > 0$
for every nonempty open subset $U \subset G^{(0)}$.

The \ca\  $C_{\mathrm{r}}^* (G_0)$ is AF by Corollary~\ref{CStarIsAF}.
It is unital,
so by Lemma~\ref{UAFHasTr} it has a normalized trace $\ta$.
Proposition~\ref{TrAndMeas} therefore implies the existence
of a $G$-invariant Borel probability measure $\mu$ on $G^{(0)}$.

Let $U \subset G^{(0)}$ be open and nonempty, and suppose that
$\mu (U) = 0$.
Let $V = r \left( s^{-1} ( U ) \right)$, using the range and
source maps of $G$.
Then $V$ is a $G$-invariant subset of $G^{(0)}$
(Definition~\ref{PDfn}(2)) which contains $U$, and it is open
because $r$ is a local \hme\  (Lemma~\ref{CGElem}(2)).
Therefore $V = G^{(0)}$.
Since every element of a Cantor groupoid $G$ is contained in a
compact open $G$-set, and since $G$ is second countable,
there exists a countable base for the topology of $G$ consisting
of compact open $G$-sets.
In particular, there is a countable collection of compact open $G$-sets,
say $W_1, \, W_2, \, \dots$, such that $s (W_n) \subset U$ for all $n$
and $V = \bigcup_{n = 1}^{\infty} r (W_n)$.
Using Lemma~\ref{GSetMeas}, we get
\[
\mu (r (W_n)) = \mu (s (W_n)) \leq \mu (U) = 0,
\]
whence
\[
\mu \ts{ \left( G^{(0)} \right) } = \mu (V)
 \leq \sum_{n = 1}^{\infty} \mu (r (W_n)) = 0.
\]
This contradicts the assumption that $\mu$ is a probability measure.
\end{pff}

\smallskip

We now show that simplicity of $C_{\mathrm{r}}^* (G)$
implies that of $C_{\mathrm{r}}^* (G_0)$.
We need a lemma.

\begin{lem}\label{AFFaithTr}
Let $A$ be a unital AF algebra.
Suppose $\ta (p) > 0$ for every normalized trace $\ta$ on $A$ and
every nonzero \pj\  $p \in A$.
Then $A$ is simple.
\end{lem}

\begin{pff}
Suppose $A$ is not simple.
Let $I$ be a nontrivial ideal in $A$.
Then $A/I$ is a unital AF algebra.
By Lemma~\ref{UAFHasTr}, there is a normalized trace $\ta$ on $A/I$.
Let $\pi : A \to A/I$ be the quotient map.
Since $I$ is AF, there is a nonzero \pj\  $p \in I$.
Then $\ta \circ \pi$ is a normalized trace on $A$ and $p$
is a nonzero \pj\  in $A$ such that $\ta \circ \pi (p) = 0$.
\end{pff}

\begin{prp}\label{AFSimple}
Let $G$ be an almost AF Cantor groupoid, with AF subgroupoid $G_0$
as in Definition~\ref{TAFGDfn}(1).
Suppose $C_{\mathrm{r}}^* (G)$ is simple.
Then $C_{\mathrm{r}}^* (G_0)$ is a simple AF algebra.
\end{prp}

\begin{pff}
The algebra $C_{\mathrm{r}}^* (G_0)$ is AF
because $G_0$ is an AF groupoid.
(See Proposition~1.15 in Chapter~3 of~\cite{Rn}.)
It is unital because the unit space of $G_0$ is compact.
By Proposition~\ref{TrAndMeas}, every trace on $C_{\mathrm{r}}^* (G_0)$
is the restriction of a trace on $C_{\mathrm{r}}^* (G)$.
Since this algebra is simple, every normalized trace on it is strictly
positive on every nonzero \pj\  in $C_{\mathrm{r}}^* (G)$, and in
particular on every nonzero \pj\  in $C_{\mathrm{r}}^* (G_0)$.
So $C_{\mathrm{r}}^* (G_0)$ is simple by Lemma~\ref{AFFaithTr}.
\end{pff}

\smallskip

We close this section with one significant unanswered question.

\begin{qst}\label{AmenQ}
Is an almost AF Cantor groupoid necessarily amenable?
\end{qst}

\begin{rmk}\label{AmenY}
If the almost AF Cantor groupoid $G$ is a \tggp\  $\Gm \times X$
(as in Example~\ref{CTGEx}), then the answer is yes; in fact,
the group $\Gm$ is necessarily amenable.
See Example~2.7(3) of~\cite{AD}.
\end{rmk}

\section{Traces and order on K-theory}\label{Sec:TOKTh}

\indent
In this section, we prove that if $G$ is an almost AF Cantor groupoid
such that $C_{\mathrm{r}}^* (G)$ is simple,
then the traces on $C_{\mathrm{r}}^* (G)$ determine the order
on $K_0 ( C_{\mathrm{r}}^* (G) )$, that is, if
$\et \in K_0 ( C_{\mathrm{r}}^* (G) )$ and $\ta_* (\et) > 0$ for all
normalized traces, then $\et > 0$.
When $G$ is the groupoid of
a substitution tiling system as in~\cite{Pt5},
this is the main result of that paper.
Theorem~\ref{PtIsAAF} below implies that such a groupoid
is in fact an almost AF Cantor groupoid.

The result of this section will be strengthened in
Section~\ref{Sec:tsr} below.

Although the proofs are a bit different (and, we hope, conceptually
simpler), the basic idea of this section is entirely contained
in Section~3 of~\cite{Pt5}.

\begin{lem}\label{LargeNorm}
Let $G$ be an almost AF Cantor groupoid, with
open AF subgroupoid $G_0 \S G$ as in Definition~\ref{TAFGDfn}(1).
Let $F \S C_{\mathrm{c}} (G)$ be a finite set, and let $\ep > 0$.
Then for every $\ep > 0$
there exists a compact open subset $V$ of $G^{(0)}$ such that, with
\[
p = \ch_V \in C \ts{\left( G^{(0)} \right)} \S C_{\mathrm{r}}^* (G_0),
\]
we have:
\begin{itemize}
\item[(1)]
$r (\supp (f) \cap [G \SM G_0])
    \cup s (\supp (f) \cap [G \SM G_0] ) \S V$
for all $f \in F$.
\item[(2)]
$\| (1 - p) f (1 - p) \| > \| f \| - \ep$ for all $f \in F$.
\item[(3)]
$\ta (p) < \ep$
for every normalized trace $\ta$ on $C_{\mathrm{r}}^* (G)$.
\end{itemize}
\end{lem}

\begin{pff}
We start by choosing $V$ so that~(1) and~(2) are satisfied.

Let $F = \{ f_1, \, f_2, \, \dots, \, f_n \}$.
By Corollary~\ref{AAFInjRepn}, there is a $G$-invariant
probability measure $\mu$ on $G^{(0)}$ whose associated regular
representation $\pi$ (see Remark~\ref{CGRmk}(3))
is faithful on $C_{\mathrm{r}}^* (G)$.
Let $\nu$ be as in Remark~\ref{CGRmk}(3).
Choose
\[
\xi_1, \, \xi_2, \, \dots, \, \xi_n,
  \, \et_1, \, \et_2, \, \dots, \, \et_n \in C_{\mathrm{c}} (G)
 \S L^2 (G, \nu)
\]
such that
\[
\| \xi_k \| = \| \et_k \| = 1 \andeqn
| \langle \pi (f_k) \xi_k, \, \et_k \rangle |
    > \| f_k \| - {\textstyle{ \frac{1}{2} }} \ep
\]
for $1 \leq k \leq n$.

Let
\[
K
= \bigcup_{k = 1}^n
  \left[ \rule{0em}{2.5ex}
      \supp (f_k) \cup \supp (\xi_k) \cup \supp (\et_k) \right].
\]
Then $K \cap (G \SM G_0)$ is a compact subset of $G \SM G_0$,
so $s (K \cap [G \SM G_0])$ is a thin set in $G^{(0)}$
by Definition~\ref{TAFGDfn}(1) and
has measure zero by Lemma~\ref{MZero}(3).
Considering $\{ g^{-1}\colon g \in K \}$ in place of $K$,
we also get $\mu (r (K \cap [G \SM G_0]) ) = 0$.
Therefore
\[
L = s (K \cap [G \SM G_0]) \cup r (K \cap [G \SM G_0] )
\]
is a compact subset of $G^{(0)}$ with $\mu (L) = 0$.
Choose a decreasing sequence of compact open sets $V_l \S G^{(0)}$
such that $\bigcap_{l = 1}^{\I} V_l = L$.
Set $p_l = \ch_{G^{(0)} \SM V_l} \in C \left( G^{(0)} \right)$.
Then one checks, for example by using the formula for
$\langle \pi (p_l) \xi, \, \et \rangle$ in Remark~\ref{CGRmk}(3),
that if $\xi \in C_{\mathrm{c}} (G) \S L^2 (G, \nu)$ then
\[
(\pi (p_l) \xi) (g) = \left\{ \begin{array}{ll}
     0       & \hspace{3em}  r (g) \in V_l \\
     f (g)   & \hspace{3em}  {\mbox{otherwise}}
    \end{array} \right..
\]
Use Lemma~\ref{GSetCover} to write $K$ as the union of compact $G$-sets
$K_1, \, K_2, \, \dots, K_N$.
{}From Lemma~\ref{GSetMeas}, and because $\mu$ is $G$-invariant, we get
\[
\nu (r^{-1} (V_l) \cap K ) = \sum_{j = 1}^N \nu (r^{-1} (V_l) \cap K_j )
  = \sum_{j = 1}^N \mu (r (r^{-1} (V_l) \cap K_j) ) \leq N \mu (V_l).
\]
Set $W_l = r^{-1} (V_l) \cap K$.
Then $W_1 \supset W_2 \supset \cdots$, we have
$\nu (W_l) \leq N \mu (V_l) \to 0$ (because $\mu (L) = 0$),
and for each $k$ we have
\[
\pi (p_l) \xi_k = \ch_{G \SM W_l} \xi_k \andeqn
\pi (p_l) \et_k = \ch_{G \SM W_l} \et_k
\]
(pointwise product on the right).
So $\pi (p_l) \xi_k \to \xi_k$ and $\pi (p_l) \et_k \to \et_k$
almost everywhere $[\nu]$ as $l \to \I$.
Applying the Dominated Convergence Theorem, we get
\[
\limi{l} \| \pi (p_l) \xi_k - \xi_k \| = 0 \andeqn
\limi{l} \| \pi (p_l) \et_k - \et_k \| = 0
\]
for $1 \leq k \leq n$.
Therefore there is $l$ such that
\[
| \langle \pi (f_k) \pi (p_l) \xi_k, \, \pi (p_l) \et_k \rangle |
    > \| f_k \| - \ep
\]
for $1 \leq k \leq n$.
So, using $\| \xi_k \| = \| \et_k \| = 1$, we get
\[
\| p_l f_k p_l \|
  \geq | \langle \pi (p_l f_k p_l) \xi_k, \, \et_k \rangle |
  = | \langle \pi (f_k) \pi (p_l) \xi_k, \, \pi (p_l) \et_k \rangle |
  > \| f_k \| - \ep.
\]
Take $V = V_l$.
With this choice, parts~(1) and~(2) of the conclusion hold.

To obtain part~(3), apply Lemma~\ref{MZero}(2) with $\ep$ as given,
and taking for $K$ the set
\[
[G \SM G_0] \cap \bigcup_{k = 1}^n [\supp (f_k) \cup \supp (f_k^*)]
\]
(which is thin in $G_0$).
Call the resulting set $W$.
Then replace $V$ by $V \cap W$.
This clearly does not affect the validity of parts~(1) and~(2),
and we get part~(3) by Proposition~\ref{TrAndMeas}.
\end{pff}

\begin{lem}\label{RR0IntPj}
Let $A$ be a \ca\  with real rank zero.
Let $a, b \in A$ be positive elements with
\[
\| a \|, \, \| b \| \leq 1 \andeqn a b = b.
\]
Let $\ep > 0$.
Then there is a \pj\  $p \in \ov{b A b}$ such that
\[
a p = p \andeqn \| p b - b \| < \ep.
\]
\end{lem}

\begin{pff}
Let $B = \ov{b A b}$.
Then $a x = x$ for all $x \in B$.
Since $A$ has real rank zero,
the \hsa\  $B$ has an approximate identity consisting of \pj s.
Since $b \in B$, there is $p \in B$ with $\| p b - b \| < \ep$.
\end{pff}

\begin{lem}\label{CPjInAF}
Let $G$ be an almost AF Cantor groupoid, with
open AF subgroupoid $G_0 \S G$ as in Definition~\ref{TAFGDfn}(1).
Let $e \in C_{\mathrm{r}}^* (G)$ be a \pj, and let $\ep > 0$.
Then there is a \pj\  $q \in C_{\mathrm{r}}^* (G_0)$ which is
\mvnt\  to a sub\pj\  of $e$ and such that $\ta (e) - \ta (q) < \ep$
for every normalized trace $\ta$ on $C_{\mathrm{r}}^* (G)$.
\end{lem}

\begin{pff}
\Wolog\  $\ep < 6$.

Choose $\dt_0 > 0$ such that whenever $A$ is a \ca\  and
$p_1, \, p_2 \in A$ are \pj s such that $\| p_1 p_2 - p_2 \| < \dt_0$,
then $p_2$ is \mvnt\  to a sub\pj\  of $p_1$.
Define a \cfn\  $f : [0, \infty) \to [0, 1]$ by
\[
f (t) = \left\{ \begin{array}{ll}
    6 \ep^{-1} t
            & \hspace{3em}  0 \leq t \leq  \ts{ \frac{1}{6} } \ep \\
    1       & \hspace{3em}  \ts{ \frac{1}{6} } \ep \leq t
    \end{array} \right..
\]
Choose $\dt > 0$ such that whenever $A$ is a \ca\  and
$a_1, \, a_2 \in A$ are positive elements with
$\| a_1 \|, \, \| a_2 \| \leq 1$ and $\| a_1 - a_2 \| < \dt$,
then $\| f ( a_1 ) - f ( a_2 ) \| < \ts{ \frac{1}{2} } \dt_0$.

Since $C_{\mathrm{c}} (G)$ is a dense *-subalgebra of
$C_{\mathrm{r}}^* (G)$, there is a selfadjoint element
$d \in C_{\mathrm{c}} (G)$ with
\[
\| e - d \|
  < \min \ts{ \left( \frac{1}{2} \dt, \, \frac{1}{6} \ep \right) }
\andeqn \| d \| \leq 1.
\]
Apply Lemma~\ref{LargeNorm} with $F = \{ d \}$, obtaining a \pj\  %
\[
p = \ch_V \in C \ts{\left( G^{(0)} \right)} \S C_{\mathrm{r}}^* (G_0)
\]
for a suitable compact open set $V \subset G^{(0)}$,
such that
\[
r (\supp (d) \cap [G \SM G_0])
    \cup s (\supp (d) \cap [G \SM G_0] ) \S V
\]
and $\ta (p) < \ep$
for every  $\ta \in T (C_{\mathrm{r}}^* (G))$, the space of
normalized traces on $C_{\mathrm{r}}^* (G)$.
Lemma~\ref{CutToAF}(3) gives
\[
(1 - p) d, \, d (1 - p) \in C_{\mathrm{r}}^* (G_0).
\]

For every $\ta \in T ( C_{\mathrm{r}}^* (G) )$, we have
\[
\ta (p e (1 - p) ) = \ta ( (1 - p) e p) = 0,
\]
so
\[
\ta ( (1 - p) e (1 - p) ) = \ta (e) - \ta (p e p)
  \geq \ta (e) - \ta (p) > \ta (e) - \ts{ \frac{1}{6} } \ep,
\]
and (using $\| d^2 - e \| < \frac{1}{3} \ep$)
\[
\ta ( d (1 - p) d ) = \ta ( (1 - p) d^2 (1 - p) )
  > \ta (e) - \ts{ \frac{1}{2} } \ep.
\]
Also, $d (1 - p) d$ is a positive element in $C_{\mathrm{r}}^* (G_0)$.

Let $f \colon [0, \infty) \to [0, 1]$ be as above, and
define \cfn s $g, h \colon [0, \infty) \to [0, 1]$ by
\[
g (t) = \left\{ \begin{array}{ll}
    0       & \hspace{1em}  0 \leq t \leq  \ts{ \frac{1}{6} } \ep \\
    6 \ep^{-1} t - 1
            & \hspace{1em}  \ts{ \frac{1}{6} } \ep
                     \leq t \leq  \ts{ \frac{1}{3} } \ep \\
    1       & \hspace{1em}  \ts{ \frac{1}{3} } \ep \leq t
    \end{array} \right.
 \,\,\,\,\,\, \andeqn  \,\,\,\,\,\,
h (t) = \left\{ \begin{array}{ll}
    t       & \hspace{1em}  0 \leq t \leq  \ts{ \frac{1}{6} } \ep \\
   \ts{ \frac{1}{6} } \ep  & \hspace{1em}  \ts{ \frac{1}{6} } \ep \leq t
    \end{array} \right..
\]
Define
\[
a = f (d (1 - p) d), \,\,\,\,\,\, b = g (d (1 - p) d),
 \andeqn c = h (d (1 - p) d).
\]
Then $a, b, c \in C_{\mathrm{r}}^* (G_0)$ are positive, and
\[
a b = b, \,\,\,\,\,\, b + c \geq d (1 - p) d, \,\,\,\,\,\,
\| a \| \leq 1, \,\,\,\,\,\, \| b \| \leq 1,
 \andeqn \| c \| \leq \ts{ \frac{1}{6} } \ep.
\]
In particular, every $\ta \in T ( C_{\mathrm{r}}^* (G) )$ satisfies
$\ta (c) \leq \ts{ \frac{1}{6} } \ep$,
so
\[
\ta (b) = \ta (b + c) - \ta (c)
   \geq \ta ( d (1 - p) d) - \ts{ \frac{1}{6} } \ep
   > \ta (e) - \ts{ \frac{2}{3} } \ep.
\]

Since $C_{\mathrm{r}}^* (G_0)$ is an AF algebra, we can apply
Lemma~\ref{RR0IntPj} to find a \pj\  %
$q \in \ov{b C_{\mathrm{r}}^* (G_0) b}$
such that $a q = q$ and $\| q b - b \| < \ts{ \frac{1}{6} } \ep$.
Then $\| q b q - b \| < \ts{ \frac{1}{3} } \ep$, so that
for every $\ta \in T ( C_{\mathrm{r}}^* (G) )$,
\[
\ta (q) \geq \ta (q b q) > \ta (b) - \ts{ \frac{1}{3} } \ep
   > \ta (e) - \ep.
\]
This is one half of the desired conclusion.

{}From $\| d - e \| < \ts{ \frac{1}{2} } \dt$ and $\| d \| \leq 1$
we get
$\| d (1 - p) d - e (1 - p) e \| < \dt$, so the choice of $\dt$ at the
beginning of the proof gives
$\| a - f ( e (1 - p) e ) \| < \ts{ \frac{1}{2} } \dt_0$.
Since
\[
e f ( e (1 - p) e ) = f ( e (1 - p) e ),
\]
we get $\| e a - a \|  < \dt_0$.
Also $a q = q$, so $\| e q - q \| < \dt_0$.
The choice of $\dt_0$ at the beginning of the proof implies that
$q$ is \mvnt\  to a sub\pj\  of $e$.
This the other half of the desired conclusion.
\end{pff}

\smallskip

The following result is a K-theoretic version of
Blackadar's Second Fundamental Comparability Question
(\cite{Bl4}, 1.3.1).
We will get the full result in Section~\ref{Sec:tsr}, after we have
proved stable rank one.

\begin{thm}\label{KFCQ2}
Let $G$ be an almost AF Cantor groupoid.
Suppose that $C_{\mathrm{r}}^* (G)$ is simple.
If $\et \in K_0 ( C_{\mathrm{r}}^* (G) )$ satisfies $\ta_* (\et) > 0$
for all normalized traces $\ta$ on $C_{\mathrm{r}}^* (G)$,
then there is a \pj\  %
\[
e \in \Mi (C_{\mathrm{r}}^* (G))
    = \bigcup_{n = 1}^{\infty} M_n (C_{\mathrm{r}}^* (G))
\]
such that $\et = [e]$.
\end{thm}

\begin{pff}
Write $\et = [q] - [p]$ for \pj s $p, q \in \Mi (C_{\mathrm{r}}^* (G))$.
Choose $n$ large enough that both $p$ and $q$ are in
$M_n (C_{\mathrm{r}}^* (G))$.
Then $G \times \{1, 2, \dots, n\}^2$, with the groupoid structure being
given by $(g, j, k) (h, k, l) = (g h, j, l)$ when $g$ and $h$ are
composable in $G$, and all other pairs not composable,
is an almost AF Cantor groupoid whose reduced \ca\  is
the simple \ca\  $M_n (C_{\mathrm{r}}^* (G))$.
Replacing $G$ by $G \times \{1, 2, \dots, n\}^2$, we may therefore
assume $p, q \in C_{\mathrm{r}}^* (G)$.

Since $C_{\mathrm{r}}^* (G)$ is simple, all traces are faithful.
Because $T ( C_{\mathrm{r}}^* (G) )$ is weak* compact, there is
$\ep > 0$ such that $\ta (q) - \ta (p) > \ep$
for all $\ta \in T ( C_{\mathrm{r}}^* (G) )$.
Apply Lemma~\ref{CPjInAF} twice, once to find a \pj\  %
$q_0 \in C_{\mathrm{r}}^* (G_0)$ which is \mvnt\  to a sub\pj\  %
$q_1$ of $q$ and such that $\ta (q) - \ta (q_0) < \frac{1}{3} \ep$ for
all $\ta \in T (C_{\mathrm{r}}^* (G) )$, and once to find a \pj\  %
$f_0 \in C_{\mathrm{r}}^* (G_0)$ which is \mvnt\  to a sub\pj\  $f_1$
of $1 - p$ and such that $\ta (1 - p) - \ta (f_0) < \frac{1}{3} \ep$
for all $\ta \in T (C_{\mathrm{r}}^* (G) )$.
Then (using Proposition~\ref{TrAndMeas}) we get
$0 < \ta (q_0) - \ta (1 - f_0) < \frac{1}{3} \ep$ for
all $\ta \in T (C_{\mathrm{r}}^* (G_0) )$.
Since $C_{\mathrm{r}}^* (G_0)$ is a simple AF algebra
(by Proposition~\ref{AFSimple}),
we can apply the last part of Corollary 6.9.2 of~\cite{Bl3}
to find a \pj\  $p_0 \leq q_0$ which is \mvnt\  to $1 - f_0$.
Now we write, in $K_0 ( C_{\mathrm{r}}^* (G) )$,
\begin{align*}
\et & = [q] - [p] = ( [q] - [q_0] ) + ( [q_0] - [p_0] ) + ([p_0] - [p])
         \\
 & = [q - q_1] + [q_0 - p_0] + [1 - p - f_1] > 0,
\end{align*}
as desired.
\end{pff}

\section{Real rank of the C*-algebra of an almost AF
  groupoid}\label{Sec:RRZ}

\indent
In this section, we prove that if $G$ is an almost AF Cantor groupoid
such that $C_{\mathrm{r}}^* (G)$ is simple,
then $C_{\mathrm{r}}^* (G)$ has real rank zero.

\begin{lem}\label{Corner}
Let $G$ be an almost AF Cantor groupoid, with
open AF subgroupoid $G_0 \S G$ as in Definition~\ref{TAFGDfn}(1).
Let $F \S C_{\mathrm{c}} (G)$ be a finite set, and let $n \in \N$.
Then there exists a compact open subset $V$ of $G^{(0)}$ such that
the \pj\  $p = \ch_V \in C \left( G^{(0)} \right)$ satisfies:
\begin{itemize}
\item[(1)]
$(1 - p) f, \, f (1 - p) \in C_{\mathrm{c}} (G_0)$
for all $f \in F$.
\item[(2)]
There are $n$ \mops\  in $C \left( G^{(0)} \right)$, each of which
is \mvnt\  to $p$ in $C_{\mathrm{r}}^* (G_0)$.
\end{itemize}
\end{lem}

\begin{pff}
Let
\[
K = \bigcup_{f \in F} [\supp (f) \cup \supp (f_k^*)].
\]
Then $K \cap (G \SM G_0)$ is a compact subset of $G \SM G_0$,
so by Definition~\ref{TAFGDfn}, the set
$s (K \cap [G \SM G_0])$ is thin in $G^{(0)}$.
By Lemma~\ref{MZero}(1), there exist a compact open set $W$
containing $K$ and compact open $G_0$-sets
$W_1, W_2, \dots, W_n \S G_0$, such that $s (W_k) = W$ and the sets
$r (W_1), \, r ( W_2), \, \dots, \, r (W_n)$ are pairwise disjoint
compact open subsets of $G^{(0)}$.
Then the functions $\ch_{W_k} \in C_{\mathrm{c}} (G_0)$ are partial
isometries which implement \mvnc s between $p$ and the $n$ \mops\  %
$\ch_{r (W_1)}, \, \ch_{r ( W_2)}, \, \dots, \, \ch_{r (W_n)}$.

We clearly have $s (\supp (f) \cap [G \SM G_0]) \S K$ for all $f \in F$.
Also,
\[
\supp (f^*) \cap [G \SM G_0]
  = \{ g^{-1} \colon g \in \supp (f) \cap [G \SM G_0] \},
\]
so
\[
r (\supp (f) \cap [G \SM G_0]) = s (\supp (f^*) \cap [G \SM G_0]) \S K
\]
for all $f \in F$.
Therefore $(1 - p) f, \, f (1 - p) \in C_{\mathrm{c}} (G_0)$
by Lemma~\ref{CutToAF}(3).
\end{pff}

\begin{lem}\label{MakePjInvFD}
Let $A$ be a finite dimensional \ca.
Let $a \in A_{\sa}$, let $p \in A$ be a \pj, and let $n \in \N$.
Then there is a \pj\  $q \in A$ such that
\[
p \leq q, \,\,\,\,\,\,\, [q] \leq 2 n [p] \in K_0 (A), \andeqn
\| q a - a q \| \leq {\textstyle{ \frac{1}{n} }} \| a \|.
\]
\end{lem}

\begin{pff}
The result is trivial if $a = 0$,
so, by scaling, we may assume $\| a \| = 1$.
\Wolog\  $A = M_m$, which we think of as operators on $\C^m$,
and $a$ is diagonal.
Making suitable replacements of the diagonal entries of $a$,
we find $b \in (M_m)_{\sa}$ such that
\[
\| a - b \| \leq {\textstyle{ \frac{1}{2n} }} \andeqn
\spec (b) \subset
\left\{ - \frac{2 n - 1}{2 n}, \,  - \frac{2 n - 3}{2 n}, \, \dots, \,
   \frac{2 n - 3}{2 n}, \, \frac{2 n - 1}{2 n} \right\}.
\]
Define $u = \exp ( \pi i b)$.
Then $u$ is a unitary in $M_m$ with $u^{2 n} = - 1$.
Moreover, with $\log$ being the continuous branch with values such that
${\mathrm{Im}} (\log (\zt)) \in (- \pi, \pi)$, we have
$b = \frac{1}{\pi i} \log (u)$.

Let $H \S \C^m$ be the linear span of the spaces $u^k p \C^m$
for $0 \leq k \leq 2 n - 1$.
Then $u H = H$ since $u^{2 n} = - 1$,
and $\dim (H) \leq 2 n \cdot \rank (p)$.
Let $q \in M_m$ be the orthogonal \pj\  onto $H$.
Then $p \leq q$ and $[q] \leq 2 n [p] \in K_0 (M_m)$.
Since $u q = q u$, functional calculus gives $b q = q b$.
Since $\| a - b \| \leq {\textstyle{ \frac{1}{2n} }}$, this gives
$\| q a - a q \| \leq {\textstyle{ \frac{1}{n} }}$.
\end{pff}

\begin{lem}\label{MakePjInvAF}
Let $A$ be a unital AF algebra.
Let $a \in A_{\sa}$ be nonzero, and let $p \in A$ be a \pj.
Let $\ep > 0$, and let $n \in \N$ satisfy $n > \frac{1}{\ep}$.
Then there is a \pj\  $q \in A$ such that
\[
p \leq q, \,\,\,\,\,\,\, [q] \leq 2 n [p] \in K_0 (A), \andeqn
\| q a - a q \| < \ep \| a \|.
\]
\end{lem}

\begin{pff}
There are finite dimensional subalgebras of $A$ which contain
\pj s arbitrarily close to $p$ and
selfadjoint elements arbitrarily close to $a$ and with the same norm.
Conjugating by suitable unitaries, we find that there are
finite dimensional subalgebras of $A$ which exactly contain $p$
and which contain
selfadjoint elements arbitrarily close to $a$ and with the same norm.
Choose such a subalgebra $B$ which contains a selfadjoint element
$b$ with
\[
\| b \| = \| a \| \andeqn
\| a - b \| < {\textstyle{ \frac{1}{2} }}
      \left( \ep - {\textstyle{ \frac{1}{n} }} \right) \| a \|.
\]
Apply Lemma~\ref{MakePjInvFD} to $B$, $b$, $p$, and $n$.
The resulting \pj\  $q \in B \S A$ has the required properties.
\end{pff}

\begin{lem}\label{FCalcApprox}
Let $r > 0$,
let $f_1, f_2, \dots, f_n : [-r, \, r] \to [0, 1]$ be \cfn s,
and let $\ep > 0$.
Then there is $\dt > 0$ such that, whenever $A$ is a unital \ca,
and whenever a \pj\  $p \in A$ and a selfadjoint element $a \in A$
satisfy $\| a \| \leq r$ and $\| p a - a p \| < \dt$, and
\[
\ta ( f_k (a)) - \ta (1 - p) > \ep
\]
for $1 \leq k \leq n$ and all normalized traces $\ta$ on $A$,
then (using functional calculus in $p A p$)
\[
\ta ( f_k (p a p)) > \ta ( f_k (a)) - \ta (1 - p) - \ep
\]
for $1 \leq k \leq n$ and all normalized traces $\ta$ on $A$.
\end{lem}

\begin{pff}
Approximating the functions $f_k$ uniformly on $[-r, r]$ by polynomials,
we can choose $\dt > 0$ so small that whenever $A$ is a unital \ca,
and whenever a \pj\  $p \in A$ and a selfadjoint element $a \in A$
satisfy $\| a \| \leq r$ and $\| p a - a p \| < \dt$, then
\[
\| f_k (p a p) - p f_k (a) p \| < \ep
\]
for $1 \leq k \leq n$.
(The expression $f_k (p a p)$ is evaluated using
functional calculus in $p A p$.)

Fix a trace $\ta$.
We estimate:
\[
\ta ( (1 - p) f_k (a) (1 - p) ) \leq \ta (1 - p) \| f_k (a) \|
   \leq \ta (1 - p)
\]
and
\[
| \ta (f_k (p a p)) - \ta ( p f_k (a) p ) |
  \leq \| f_k (p a p) - p f_k (a) p \| < \ep
\]
for $1 \leq k \leq n$.
Also,
\[
\ta ( (1 - p) f_k (a) p ) = \ta ( p f_k (a) (1 - p) ) = 0,
\]
because $\ta$ is a trace.
So
\begin{align*}
\ta ( f_k (p a p))
   & > \ta ( p f_k (a) p ) - \ep  \\
   & = \ta ( f_k (a) ) - \ep \\
   & \hspace{2.5em}
        - \left[ \rule{0em}{2.5ex}
          \ta ( (1 - p) f_k (a) (1 - p) ) + \ta ( (1 - p) f_k (a) p )
        + \ta ( p f_k (a) (1 - p) ) \right]  \\
   & \geq \ta ( f_k (a)) - \ta (1 - p) - \ep,
\end{align*}
as desired.
\end{pff}

\begin{lem}\label{DiagConj}
Let $X \S \R$ be a compact set and let $\ep > 0$.
Let $\ld_1 < \ld_2 < \cdots < \ld_n$ be real numbers such that
for every $\ld \in X$ there is $k$ with $| \ld - \ld_k | < \ep$.
Define $a \in C (X, M_{n + 1})$ by
\[
a (\ld) = \diag (\ld, \, \ld_1, \, \ld_2, \, \dots, \, \ld_n )
\]
for $\ld \in X$.
Then there exists a selfadjoint diagonal element
$b \in C (X, M_{n + 1})$
and a unitary $u \in C (X, M_{n + 1})$ such that
\[
\spec (b) \S \{ \ld_1, \, \ld_2, \, \dots, \, \ld_n \}
\andeqn \| u a u^* - b \| < 2 \ep.
\]
\end{lem}

\begin{pff}
We first prove this in the case that $X$ is an interval $[\af, \bt]$
and $\ld_{k + 1} - \ld_k < 2 \ep$ for $1 \leq k \leq n - 1$.
Let
\[
\rh = \sup_{\ld \in [\af, \bt]}
   \dist (\ld, \, \{ \ld_1, \, \ld_2, \, \dots, \, \ld_n \} ).
\]
Then $\rh < \ep$ by assumption.
Choose numbers $\ep_k$ with $0 < \ep_k < \ep - \rh$
such that
\[
\ld_1 + \ep_1 < \ld_2 - \ep_2 < \ld_2 + \ep_2 < \cdots
  < \ld_{n - 1} + \ep_{n - 1} < \ld_n - \ep_n.
\]
For any $k$, let $I_k$ denote the $k \times k$ identity matrix.
For $0 \leq k \leq n$, define a unitary $v_k \in M_{n + 1}$ by
\[
v_k = \left( \begin{array}{ccccccc}
0       & 1       & 0       & \cdots  & \cdots  & 0       & 0        \\
0       & 0       & 1       & \cdots  & \cdots  & 0       & 0        \\
\vdots  & \vdots  & \vdots  & \ddots  &         & \vdots  & \vdots   \\
\vdots  & \vdots  & \vdots  &         & \ddots  & \vdots  & \vdots   \\
0       & 0       & 0       & \cdots  & \cdots  & 1       & 0        \\
0       & 0       & 0       & \cdots  & \cdots  & 0       & 1        \\
1       & 0       & 0       & \cdots  & \cdots  & 0       & 0
\end{array} \right)
\oplus I_{n - k},
\]
where the first block is a cyclic shift
of size $(k + 1) \times (k + 1)$.
Note that
\[
v_k \diag (\ld, \, \ld_1, \, \ld_2, \, \dots, \, \ld_n ) v_k^*
  = \diag (\ld_1, \, \dots, \, \ld_k, \, \ld,
               \, \ld_{k + 1}, \, \dots, \, \ld_n ).
\]
Choose a \ct\  path $\ld \mapsto z_k (\ld)$
in the unitary group of $M_2$,
for $\ld_k - \ep_k \leq \ld \leq \ld_k + \ep_k$, such that
\[
z_k (\ld_k - \ep_k) = 1 \andeqn
 z_k (\ld_k + \ep_k)
  = \left( \begin{array}{cc} 0 & 1 \\ 1 & 0 \end{array} \right).
\]
Define
\[
w_k ( \ld) = I_{k - 1} \oplus z_k (\ld) \oplus I_{n - k}
   \in M_{n + 1}
\]
for $\ld_k - \ep_k \leq \ld \leq \ld_k + \ep_k$.
Define
\[
u (\ld) = \left\{ \begin{array}{ll}
     v_1     & \hspace{3em}  \ld \leq \ld_1 + \ep_1 \\
     v_k     & \hspace{3em}
          \ld_k + \ep_k \leq \ld \leq \ld_{k + 1} - \ep_{k + 1}
              \,\, {\mbox{with}} \,\, 1 \leq k \leq n - 1 \\
     w_k (\ld) v_{k - 1}  & \hspace{3em}
          \ld_k - \ep_k \leq \ld \leq \ld_k + \ep_k
              \,\, {\mbox{with}} \,\, 2 \leq k \leq n \\
     v_n     & \hspace{3em}  \ld_n + \ep_n \leq \ld
         \end{array} \right..
\]

Set
\[
b_0 = \diag (\ld_1, \, \ld_1, \, \ld_2, \, \dots, \, \ld_n ).
\]
We claim that
\[
\| u (\ld) a (\ld) u (\ld)^* - b_0 \| < 2 \ep
\]
for all $\ld \in [\af, \bt]$.
This will prove the special case, with $b$ being the constant
function $b (\ld) = b_0$ for all $\ld \in [\af, \bt]$.
For $1 \leq k \leq n - 1$ and
$\ld_k + \ep_k \leq \ld \leq \ld_{k + 1} - \ep_{k + 1}$, this
difference is
\begin{align*}
& \| \diag (\ld_1, \, \ld_2, \, \ld_3, \, \dots, \, \ld_k, \, \ld,
               \, \ld_{k + 1}, \, \dots, \, \ld_n )   \\
& \hspace*{10em}
      - \diag (\ld_1, \, \ld_1, \, \ld_2, \, \dots, \, \ld_{k - 1},
               \, \ld_k, \, \ld_{k + 1}, \, \dots, \, \ld_n ) \|,
\end{align*}
which is at most
\[
\max_{1 \leq j \leq k - 1} (\ld_{j + 1} - \ld_j ) < 2 \ep.
\]
Similar calculations show that
$\| u (\ld) a (\ld) u (\ld)^* - b_0 \| < 2 \ep$
for $\ld \leq \ld_1 + \ep_1$ and for $\ld_n + \ep_n \leq \ld$.
Using the same method for $2 \leq k \leq n$ and
$\ld_k - \ep_k \leq \ld \leq \ld_k + \ep_k$, we get
\begin{align*}
& \| u (\ld) a (\ld) u (\ld)^* - b_0 \|  \\
& \hspace*{1em}
  \leq \max \left( \max_{1 \leq j \leq k - 2} (\ld_{j + 1} - \ld_j ),
    \,\, \left\| z_k (\ld)
       \left( \begin{array}{cc} \ld & 0 \\ 0 & \ld_k \end{array} \right)
        z_k (\ld)^*
     - \left( \begin{array}{cc}
          \ld_{k - 1} & 0 \\ 0 & \ld_k \end{array} \right)
        \right\| \right).
\end{align*}
We have already seen that the first term is less than $2 \ep$.
To estimate the other term, first observe that
\begin{align*}
\left\| z_k (\ld)
       \left( \begin{array}{cc} \ld & 0 \\ 0 & \ld_k \end{array} \right)
        z_k (\ld)^*
   - \left( \begin{array}{cc} \ld_k & 0 \\ 0 & \ld_k \end{array} \right)
        \right\|
& = \left\|
       \left( \begin{array}{cc} \ld & 0 \\ 0 & \ld_k \end{array} \right)
   - \left( \begin{array}{cc} \ld_k & 0 \\ 0 & \ld_k \end{array} \right)
        \right\|  \\
&  < \ep - \rh
\end{align*}
because $z_k (\ld)$ commutes with
$\left( \begin{array}{cc} \ld_k & 0 \\ 0 & \ld_k \end{array} \right)$.
Therefore this term
is less than
\[
(\ep - \rh) + (\ld_k - \ld_{k - 1}) \leq \ep - \rh + 2 \rh < 2 \ep.
\]
So $\| u (\ld) a (\ld) u (\ld)^* - b_0 \| < 2 \ep$
for all $\ld \in [\af, \bt]$.
The special case is now proved.

In the general case, set
\[
\rh = \sup_{\ld \in X}
   \dist (\ld, \, \{ \ld_1, \, \ld_2, \, \dots, \, \ld_n \} ) < \ep
\]
and
\[
Y = \{ \ld \in \R \colon
 \dist (\ld, \, \{ \ld_1, \, \ld_2, \, \dots, \, \ld_n \} ) \leq \rh \}.
\]
Then $X \S Y$, and it suffices to prove the result for
$Y$ in place of $X$.
Now $Y$ is the finite disjoint union of compact intervals.
For each such interval $J$, there are $k \leq l$ such that
\begin{itemize}
\item[(1)]
$\{ \ld_k, \, \dots, \, \ld_l \}
   = J \cap \{ \ld_1, \, \dots, \, \ld_n \}$.
\item[(2)]
$\ld_{j + 1} - \ld_j < 2 \ep$ for $k \leq j \leq l - 1$.
\item[(3)]
For every $\ld \in J$ there is $j$ such that $k \leq j \leq l$
and $| \ld - \ld_k | < \ep$.
\end{itemize}
Apply the result of the special case, obtaining $a_J$, $b_J$, and $u_J$.
Then on $J$ define
\[
a_0 (\ld) = \diag (\ld_1, \, \dots, \, \ld_{k - 1} )
  \oplus a_J (\ld)
  \oplus \diag (\ld_{l + 1}, \, \dots, \, \ld_n ),
\]
\[
b_0 (\ld) = \diag (\ld_1, \, \dots, \, \ld_{k - 1} )
  \oplus b_J (\ld)
  \oplus \diag (\ld_{l + 1}, \, \dots, \, \ld_n ),
\]
and
\[
u_0 (\ld) = I_{k - 1} \oplus u_J (\ld) \oplus I_{n - l}.
\]
A suitable permutation matrix $w$ conjugates $a_0 (\ld)$ to $a (\ld)$
as defined in the statement,
and we set $u = w u_0 w^*$ and $b = w b_0 w^*$.
\end{pff}

\smallskip

The following lemma gives a suitable version of a now standard
technique, which goes back at least to the proof of
Lemma~1.7 of~\cite{Cu2}.

\begin{lem}\label{ApproxCorner}
Let $A$ be a \ca, let $a \in A$ be normal, let $\ld_0 \in \C$,
let $\ep > 0$, and let $f : \C \to [0, 1]$ be a \cfn\  such that
$\supp (f)$ is contained in the open $\ep$-ball with center $\ld_0$.
Let $p$ be a \pj\  in the \hsa\  generated by $f (a)$.
Then
\[
\| p a - a p \| < 2 \ep \andeqn \| p a p - \ld_0 p \| < \ep.
\]
\end{lem}

\begin{pff}
Choose a \cfn\  $g : \C \to \C$ such that $| g ( \ld) - \ld | < \ep$
for all $\ld \in \C$ and $g (\ld) = \ld_0$ for all $\ld \in \supp (f)$.
Then
\[
\| g (a) - a \| < \ep, \,\,\,\,\,\, p g (a) = g (a) p, \andeqn
p g (a) p = \ld_0 p.
\]
\end{pff}

\begin{thm}\label{RR0}
Let $G$ be an almost AF Cantor groupoid.
Suppose that $C_{\mathrm{r}}^* (G)$ is simple.
Then $C_{\mathrm{r}}^* (G)$ has real rank zero
in the sense of~\cite{BP}.
\end{thm}

\begin{pff}
Let $a \in C_{\mathrm{r}}^* (G)$ be selfadjoint with $\| a \| \leq 1$.
We approximate $a$ by a selfadjoint element with finite spectrum.
Since $C_{\mathrm{c}} (G)$ is a dense *-subalgebra of
$C_{\mathrm{r}}^* (G)$, we can approximate $a$ arbitrarily well by an
element $a_0 \in C_{\mathrm{c}} (G)$; replacing $a_0$ by
$\frac{1}{2} (a_0 + a_0^*)$ and suitably scaling, we
can approximate $a$ arbitrarily well by a selfadjoint
element $a_0 \in C_{\mathrm{c}} (G)$ with $\| a_0 \| \leq 1$.
Thus, \wolog\  $a \in C_{\mathrm{c}} (G)$.

Let $\ep > 0$.
Choose real numbers $\ld_1 < \ld_2 < \cdots < \ld_n$ in $\spec (a)$
such that for every $\ld \in \spec (a)$ there is $k$ such that
$| \ld - \ld_k | < \frac{1}{4} \ep$.
Set
\[
\ep_0 = \frac{\ep}{2 (2 n^2 + 2 n + 3)},
\]
and choose \cfn s $f_k, \, g_k : [-1, 1] \to [0, 1]$
such that the $f_k$ have disjoint supports and
\[
g_k (\ld_k) = 1, \,\,\,\,\,\, f_k g_k = g_k,
\andeqn
\supp (f_k) \S \{ \ld \in \R \colon | \ld - \ld_k | < \ep_0 \}
\]
for $1 \leq k \leq n$.

Let $T (C_{\mathrm{r}}^* (G) )$ be the set of normalized traces
on $C_{\mathrm{r}}^* (G)$, and define
\[
\af = \inf_{1 \leq k \leq n}
  \left( \inf_{\ta \in T (C_{\mathrm{r}}^* (G) )}
    \ta (g_k (a)) \right).
\]
Each $g_k (a)$ is a nonzero positive element,
all traces are faithful because $C_{\mathrm{r}}^* (G)$ is simple,
and $T (C_{\mathrm{r}}^* (G) )$ is weak* compact.
Therefore $\af > 0$.

Choose $\dt > 0$ as in Lemma~\ref{FCalcApprox}, with
$\frac{1}{4} \af$ in place of $\ep$, with $r = 2$, and with
$g_1, g_2, \dots, g_n$ in place of $f_1, f_2, \dots, f_n$.
We also require $\dt < \ep_0$.
Choose $m \in \N$ with $\frac{2}{m} < \dt$, and use Lemma~\ref{Corner}
to find a \pj\  $p_0 \in C \left( G^{(0)} \right)$ which is \mvnt\  in
$C_{\mathrm{r}}^* (G_0)$ to more than $8 m \af^{-1}$ \mops\  in
$C \left( G^{(0)} \right)$, and such that
$(1 - p_0) a, \, a (1 - p_0) \in C_{\mathrm{c}} (G_0)$.
In particular,
\[
\ta (p_0) < \frac{\af}{8 m}
\]
for all $\ta \in T (C_{\mathrm{r}}^* (G_0) )$.

Define $b = a - p_0 a p_0$, which is a selfadjoint element of
$C_{\mathrm{c}} (G_0)$ with $\| b \| \leq 2$.
Because $\frac{2}{m} < \dt$, we can apply
Lemma~\ref{MakePjInvAF} with $\frac{1}{2} \dt$ in place of $\ep$
to obtain a \pj\  $p \in C_{\mathrm{r}}^* (G_0)$
such that $\| p b - b p \| < \dt$, $p_0 \leq p$,
and $[p] \leq 2 m [p_0]$ in $K_0 ( C_{\mathrm{r}}^* (G) )$.
Now $p$ commutes with $b - a = p_0 a p_0$, so also
$\| p a - a p \| < \dt$.
Furthermore, because $p \in C_{\mathrm{r}}^* (G_0)$ and $p \geq p_0$,
we get $(1 - p) a, \, a (1 - p) \in C_{\mathrm{r}}^* (G_0)$.

Define $a_0 = (1 - p) a (1 - p)$.
For every $\ta \in T (C_{\mathrm{r}}^* (G_0) )$, we have
\[
\ta (p) \leq 2 m \ta (p_0) < 2 m \cdot \frac{\af}{8 m}
 = {\textstyle{\frac{1}{4}}} \af.
\]
By the choice of $\dt$ and using Lemma~\ref{FCalcApprox}
(with $1 - p$ in place of $p$ and ${\textstyle{\frac{1}{4}}} \af$
in place of $\ep$), we get
\[
\ta (g_k ( a_0 )) >
    \ta ( g_k (a)) - \ta (p) - {\textstyle{\frac{1}{4}}} \af
  \geq
    \af - {\textstyle{\frac{1}{4}}} \af - {\textstyle{\frac{1}{4}}} \af
  = {\textstyle{\frac{1}{2}}} \af
\]
for $1 \leq k \leq n$ and all $\ta \in T (C_{\mathrm{r}}^* (G_0) )$.
Also $f_k (a_0) g_k (a_0) = g_k (a_0)$
for all $k$, and $C_{\mathrm{r}}^* (G_0)$ is an AF algebra,
so Lemma~\ref{RR0IntPj} provides \pj s $q_k \in C_{\mathrm{r}}^* (G_0)$
such that
\[
q_k \in \ov{g_k (a_0) A g_k (a_0)}, \,\,\,\,\,\,\,
f_k (a_0) q_k = q_k,
\andeqn \| q_k g_k (a_0) - g_k (a_0) \| < {\textstyle{\frac{1}{8}}} \af.
\]
Therefore
\[
\| q_k g_k (a_0) q_k - g_k (a_0) \| < {\textstyle{\frac{1}{4}}} \af.
\]

For all $\ta \in T (C_{\mathrm{r}}^* (G_0) )$, we have
\[
\ta (q_k g_k (a_0) q_k) \leq \ta (q_k)
\]
because $\| g_k (a_0) \| \leq 1$.
Combining this with the previous estimates, it follows that
\[
\ta (q_k) > \ta ( g_k (a_0)) - {\textstyle{\frac{1}{4}}} \af
  > {\textstyle{\frac{1}{2}}} \af - {\textstyle{\frac{1}{4}}} \af
  = {\textstyle{\frac{1}{4}}} \af.
\]
Since $C_{\mathrm{r}}^* (G_0)$ is simple
(by Proposition~\ref{AFSimple}),
since traces determine order on the $K_0$ group of a simple AF
algebra (see the last part of Corollary 6.9.2 of~\cite{Bl3}),
and since $\ta (p) < {\textstyle{\frac{1}{4}}} \af$ for all
$\ta \in T (C_{\mathrm{r}}^* (G_0) )$,
it follows that there are \pj s $p_k \leq q_k$
with $p_k$ \mvnt\  to $p$.
In particular, $p_k$ is in the \hsa\  of $C_{\mathrm{r}}^* (G_0)$
generated by $g_k (a_0)$.
Since the $g_k$ have disjoint supports, the $p_k$ are orthogonal;
since $a_0$ is orthogonal to $p$, so are the $p_k$.
Moreover, Lemma~\ref{ApproxCorner} implies that
\[
\| p_k a_0 - a_0 p_k \| < 2 \ep_0
\andeqn \| p_k a_0 p_k - \ld_k p_k \| < \ep_0.
\]

Define
\[
e = 1 - p - \sum_{k = 1}^n p_k.
\]
We estimate
\[
\left\| a - \left(e a e + p a p + \sum_{k = 1}^n \ld_k p_k \right)
    \right\|.
\]
First,
\[
a - a_0 - p a p = p a (1 - p) + (1 - p) a p.
\]
Since $\| p a - a p \| < \dt \leq \ep_0$,
we get
\[
\| a - a_0 - p a p \| < 2 \ep_0.
\]
Next,
\[
a_0 - \left(e a e + \sum_{k = 1}^n p_k a p_k \right)
 = \sum_{k = 1}^n p_k a_0 e + \sum_{k = 1}^n e a_0 p_k
    + \sum_{j \neq k} p_j a_0 p_k.
\]
Since $\| p_k a_0 - a_0 p_k \| < 2 \ep_0$,
each of the terms has norm at most $2 \ep_0$.
There are $2 n + (n^2 - n) = n^2 + n$ terms, so
\[
\left\| a - \left(e a e + p a p + \sum_{k = 1}^n p_k a p_k \right)
    \right\| < [ 2 (n^2 + n) + 2] \ep_0.
\]
Moreover, $\| p_k a_0 p_k - \ld_k p_k \| < \ep_0$ and the terms
$p_k a_0 p_k - \ld_k p_k$ are orthogonal, so it follows that
\[
\left\| a - \left(e a e + p a p + \sum_{k = 1}^n \ld_k p_k \right)
    \right\| < [ 2 (n^2 + n) + 3] \ep_0 = {\textstyle{\frac{1}{2}}} \ep.
\]

There is a C* subalgebra $B$ of $(1 - e) C_{\mathrm{r}}^* (G) (1 - e)$
isomorphic to $C (\spec (p a p), M_{n + 1})$ with identity
\[
1 - e = p - \sum_{k = 1}^n p_k
\]
and containing
\[
p a p + \sum_{k = 1}^n \ld_k p_k.
\]
Using Lemma~\ref{DiagConj} and the choice of the numbers $\ld_k$,
we can construct a selfadjoint element
$b \in (1 - e) C_{\mathrm{r}}^* (G) (1 - e)$ with finite spectrum
such that
\[
\left\| b - \left( p a p + \sum_{k = 1}^n \ld_k p_k \right) \right\|
 < {\textstyle{\frac{1}{2}}} \ep.
\]
Also, $e a e \in e C_{\mathrm{r}}^* (G_0) e$, which is an AF algebra,
so there is a selfadjoint element
$c \in e C_{\mathrm{r}}^* (G_0) e$ with finite spectrum such that
\[
\left\| c - e a e \right\| < {\textstyle{\frac{1}{2}}} \ep.
\]
It follows that $b + c$ is a selfadjoint element
in $C_{\mathrm{r}}^* (G)$ with finite spectrum
such that
\[
\left\| a - (b + c) \right\| < \ep.
\]
This completes the proof.
\end{pff}

\section{Stable rank of the C*-algebra of an almost AF
groupoid}\label{Sec:tsr}

\indent
In this section, we prove that if $G$ is an almost AF Cantor groupoid
such that $C_{\mathrm{r}}^* (G)$ is simple,
then $C_{\mathrm{r}}^* (G)$ has stable rank one.
This implies that the \pj s in $M_{\infty} (C_{\mathrm{r}}^* (G) )$
satisfy cancellation, and allows us to strengthen the conclusion of
Theorem~\ref{KFCQ2} to the full version of
Blackadar's Second Fundamental Comparability Question.

\begin{lem}\label{ZDOrth}
Let $G$ be an almost AF Cantor groupoid, with
open AF subgroupoid $G_0 \S G$ as in Definition~\ref{TAFGDfn}(1).
For every $\ep > 0$ there is $\dt > 0$ such that if
$a \in C_{\mathrm{c}} (G)$ satisfies $\| a \| \leq 1$ and if
$e_0 \in C_{\mathrm{r}}^* (G)$ is a nonzero \pj\  such that
$\| a e_0 \| < \dt$, then there exists a nonzero \pj\  %
$e \in C_{\mathrm{r}}^* (G_0)$ and a compact open subset
$V \S G^{(0)}$ such that,
with $p = \ch_V \in C_{\mathrm{r}}^* (G_0)$, we have:
\begin{itemize}
\item[(1)]
$r (\supp (a) \cap [G \SM G_0])
    \cup s (\supp (a) \cap [G \SM G_0] ) \S V$.
\item[(2)]
$\| a e \| < \ep$.
\item[(3)]
$e$ and $p$ are orthogonal.
\end{itemize}
\end{lem}

\begin{pff}
Choose $\dt = \min \left( \frac{1}{3} \ep, \, \frac{1}{2} \right) > 0$.
Let $a \in C_{\mathrm{c}} (G)$ satisfy $\| a \| \leq 1$ and let
$e_0 \in C_{\mathrm{r}}^* (G)$ be a nonzero \pj\  such that
$\| a e_0 \| < \dt$.
Choose $c \in C_{\mathrm{c}} (G)$ with $\| c \| = 1$ and
$\| c - e_0 \| < \dt^2$.
Apply Lemma~\ref{LargeNorm} with
$F = \{ a, \, c, \, c^*, \, c^* c \}$,
and with $\dt^2$ in place of $\ep$.
We obtain a compact open subset $V \S G^{(0)}$ such that,
with $p = \ch_V$,
Condition~(1) is satisfied, and also
\[
\| (1 - p) c^* c (1 - p) \| > 1 - \dt^2.
\]
Moreover, Lemma~\ref{CutToAF}(3)
implies that $(1 - p) c^*, \, c (1 - p) \in C_{\mathrm{r}}^* (G_0)$.
It follows that $(1 - p) c^* c (1 - p) \in C_{\mathrm{r}}^* (G_0)$,
which is an AF algebra.

Choose a \cfn\  $f$ such that
\[
[1 - \dt^2, \, 1] \S \supp (f) \S (1 - 2 \dt^2, \, 1 + \dt^2).
\]
Choose a nonzero \pj\  $e$ in the \hsa\  generated by
\[
f ( (1 - p) c^* c (1 - p) )
\]
in the AF algebra $C_{\mathrm{r}}^* (G_0)$.
Apply Lemma~\ref{ApproxCorner} in $C_{\mathrm{r}}^* (G_0)$
with this $f$ and with $\ld_0 = 1$, getting
\[
\| e - e (1 - p) c^* c (1 - p) e \| < 2 \dt^2.
\]
Since $0 \not\in \supp (f)$, we have
$e \in (1 - p) C_{\mathrm{r}}^* (G_0) (1 - p)$.
So $p$ and $e$ are orthogonal, which is Part~(3).
We combine this with the estimate above and the estimate
\[
\| c^* c - e_0 \|
   \leq \| c^* \| \| c - e_0 \| + \| c^* - e_0 \| \| e_0 \| < 2 \dt^2
\]
to obtain
\[
\| e - e_0 e \|^2 = \| e ( 1 - e_0 ) e \|
  < \| e (1 - c^* c ) e \| + 2 \dt^2 = \| e - e c^* c e \| + 2 \dt^2
  < 4 \dt^2.
\]
Therefore, since $\| a \| \leq 1$,
\[
\| a e \| \leq \| a \| \| e - e_0 e \| + \| a e_0 \| \| e \|
  < 2 \cdot \dt + \dt \leq \ep.
\]
This is Part~(2).
\end{pff}

\begin{thm}\label{tsr1}
Let $G$ be an almost AF Cantor groupoid.
Suppose that $C_{\mathrm{r}}^* (G)$ is simple.
Then $C_{\mathrm{r}}^* (G)$ has
(topological) stable rank one in the sense of~\cite{Rf}.
\end{thm}

\begin{pff}
We are going to show that every two sided zero divisor in
$C_{\mathrm{r}}^* (G)$ is a limit of invertible elements.
That is, if $a \in C_{\mathrm{r}}^* (G)$ and there are nonzero
$x, y \in C_{\mathrm{r}}^* (G)$ such that $x a = a y = 0$,
then we show that for every $\ep > 0$ there is an invertible element
$c \in C_{\mathrm{r}}^* (G)$ such that $\| a - c \| < \ep$.
It will follow from Theorem~3.3(a) of~\cite{Rr0}
(see Definition~3.1 of~\cite{Rr0})
that any $a \in C_{\mathrm{r}}^* (G)$ which is not a limit of
invertible elements is left or right invertible but not both.
Since $C_{\mathrm{r}}^* (G)$ is simple,
the definition of an almost AF Cantor groupoid implies it has a
faithful trace, and there are no such elements.

So let $a \in C_{\mathrm{r}}^* (G)$,
let $x, y \in C_{\mathrm{r}}^* (G)$ be nonzero elements
such that $x a = a y = 0$, and let $\ep > 0$.
\Wolog\  $\| a \| \leq 1$.
Since $C_{\mathrm{r}}^* (G)$ has real rank zero (Theorem~\ref{RR0}),
there are nonzero \pj s
\[
e_0 \in \ov{x^* C_{\mathrm{r}}^* (G) x} \andeqn
f_0 \in \ov{y C_{\mathrm{r}}^* (G) y^*},
\]
and we have $e_0 a = a f_0 = 0$.
Choose $\dt > 0$ in Lemma~\ref{ZDOrth} for $\frac{1}{4} \ep$ in place
of $\ep$.
Choose $b \in C_{\mathrm{c}} (G)$ such that $\| b \| \leq 1$ and
$\| a - b \| < \min \left( \dt, \, \frac{1}{4} \ep \right)$.
Then $\| e_0 b \|, \, \| b f_0 \| < \dt$.
Let
\[
K = r (\supp (b) \cap [G \SM G_0])
    \cup s (\supp (b) \cap [G \SM G_0] ).
\]
Apply Lemma~\ref{ZDOrth} to $b^*$ and $e_0$, obtaining a nonzero \pj\  %
$e_1 \in C_{\mathrm{r}}^* (G_0)$ and a compact open subset
$V \S G^{(0)}$ such that,
with $g = \ch_V \in C_{\mathrm{r}}^* (G)$, we have:
\[
K \S V, \,\,\,\,\,\,
\| e_1 b \| < \ts{ \frac{1}{4} } \ep, \andeqn e_1 g = 0.
\]
Apply Lemma~\ref{ZDOrth} to $b$ and $f_0$, obtaining a nonzero \pj\  %
$f_1 \in C_{\mathrm{r}}^* (G_0)$ and a compact open subset
$W \S G^{(0)}$ such that,
with $h = \ch_W \in C_{\mathrm{r}}^* (G)$, we have:
\[
K \S W, \,\,\,\,\,\,
\| b f_1 \| < \ts{ \frac{1}{4} } \ep, \andeqn f_1 h = 0.
\]

Choose $\rh > 0$ with
\[
\rh < \min \left(
  \inf_{\ta \in T ( C_{\mathrm{r}}^* (G) ) } \ta (e_1), \,\,
  \inf_{\ta \in T ( C_{\mathrm{r}}^* (G) ) } \ta (f_1) \right).
\]
(As usual,
$T ( C_{\mathrm{r}}^* (G) )$ is the space of normalized traces.)
The set $K$ is thin in $G_0$,
so Lemma~\ref{MZero}(2) provides a compact
open set $Z \S G^{(0)}$ such that $K \S Z$ and $\mu (Z) < \rh$
for every $G_0$-invariant Borel probability measure $\mu$ on $G^{(0)}$.
Let $p = \ch_{V \cap W \cap Z} \in C \left( G^{(0)} \right)$.
Then $p$, $e_1$, and $f_1$ are all in $C_{\mathrm{r}}^* (G_0)$,
and $e_1 p = f_1 p = 0$.
Proposition~\ref{TrAndMeas} implies that
\[
\ta (p) < \ta (e_1) \andeqn \ta (p) < \ta (f_1)
\]
for all $\ta \in T (C_{\mathrm{r}}^* (G_0))$.
Since $C_{\mathrm{r}}^* (G_0)$ is a simple AF algebra
(by Proposition~\ref{AFSimple}), it follows
(see the last part of Corollary 6.9.2 of~\cite{Bl3})
that there are \pj s $e_2, f_2 \in C_{\mathrm{r}}^* (G_0)$
which are unitarily equivalent to $p$ and with
$e_2 \leq e_1$ and $f_2 \leq f_1$.
Furthermore, since $e_2$ and $f_2$ are both orthogonal to $p$,
and $(1 - p) C_{\mathrm{r}}^* (G_0) (1 - p)$ is a simple AF algebra,
there is a unitary $w \in C_{\mathrm{r}}^* (G_0)$ such that
\[
w e_2 w^* = f_2 \andeqn w p w^* = p;
\]
also, it is easy to find a unitary $v \in C_{\mathrm{r}}^* (G_0)$
such that
\[
v f_2 v^* = p \andeqn v p v^* = f_2.
\]
Set $u = v w$, which is a unitary in $C_{\mathrm{r}}^* (G_0)$ such that
\[
u e_2 u^* = p \andeqn u p u^* = f_2.
\]

Define
\[
b_0 = (1 - e_2) b (1 - f_2).
\]
Since
\[
b_0 = b - e_2 b (1 - f_2) - b f_2,
\]
and since
\[
\| e_2 b (1 - f_2) \| \leq \| e_1 b \| < \ts{ \frac{1}{4} } \ep \andeqn
\| b f_2 \| \leq \| b f_1 \| < \ts{ \frac{1}{4} } \ep,
\]
we get $\| b - b_0 \| < \ts{ \frac{1}{2} } \ep$,
so $\| a - b_0 \| < \ts{ \frac{3}{4} } \ep$.
We have $(1 - p) b (1 - p) \in C_{\mathrm{r}}^* (G_0)$ by
Lemma~\ref{CutToAF}(3).
With respect to the decomposition of the identity
\[
1 = p + [1 - p - f_2] + f_2,
\]
we claim that $u b_0$ has the block matrix form
\[
u b_0 = \left( \begin{array}{ccc}
     0       &  0       &  0        \\
     x       &  d_0     &  0        \\
     y       &  z       &  0
    \end{array} \right)
\]
with
\[
d_0 \in (1 - p - f_2) C_{\mathrm{r}}^* (G_0) (1 - p - f_2).
\]
To see this, we observe that $u b_0 f_2 = 0$ because $b_0 f_2 = 0$,
that
\[
p u b_0 = u (u^* p u) b_0 = u e_2 b_0 = 0
\]
because $e_2 b_0 = 0$, and that
\begin{align*}
d_0 & = (1 - p - f_2) u b_0 (1 - p - f_2)
    = u (1 - e_2 - p) b_0 (1 - p - f_2) \\
  & = u [1 - e_2] [(1 - p) b_0 (1 - p)] [1 - f_2],
\end{align*}
which is in $C_{\mathrm{r}}^* (G_0)$ because $u$ and each of the
three terms in brackets are in $C_{\mathrm{r}}^* (G_0)$.
Since $C_{\mathrm{r}}^* (G_0)$ is an AF algebra, there exists
an invertible element $d \in C_{\mathrm{r}}^* (G_0)$ such that
$\| d - d_0 \| < \ts{ \frac{1}{4} } \ep$.
Then
\[
c = u^* \cdot \left( \begin{array}{ccc}
\ts{ \frac{1}{4} } \ep &  0       &  0        \\
               x       &  d       &  0        \\
               y       &  z       & \ts{ \frac{1}{4} } \ep
    \end{array} \right)
\]
is an invertible element of $C_{\mathrm{r}}^* (G)$ such that
$\| b_0 - c \| \leq \ts{ \frac{1}{4} } \ep$, so
$\| a - c \| < \ep$.
\end{pff}

\smallskip

Recall that $\Mi (A) = \bigcup_{n = 1}^{\infty} M_n (A)$.

\begin{cor}\label{Canc}
Let $G$ be an almost AF Cantor groupoid.
Suppose that $C_{\mathrm{r}}^* (G)$ is simple.
Then the \pj s in $\Mi (C_{\mathrm{r}}^* (G))$ satisfy cancellation:
if $e, f, p \in \Mi (C_{\mathrm{r}}^* (G))$ are \pj s such that
$e \oplus p$ is \mvnt\  to $f \oplus p$, then $e$ is \mvnt\  to $f$.
\end{cor}

\begin{pff}
This follows from the fact that $C_{\mathrm{r}}^* (G)$ has
stable rank one (Theorem~\ref{tsr1}) and real rank zero
(Theorem~\ref{RR0}), using Proposition 6.5.1 of~\cite{Bl3}.
\end{pff}

\begin{cor}\label{FCQ2}
Let $G$ be an almost AF Cantor groupoid.
Suppose that $C_{\mathrm{r}}^* (G)$ is simple.
Let $p, q \in \Mi (C_{\mathrm{r}}^* (G))$
be \pj s such that $\ta (p) < \ta (q)$
for all normalized traces $\ta$ on $C_{\mathrm{r}}^* (G)$.
Then $p$ is \mvnt\  to a sub\pj\  of $q$.
\end{cor}

\begin{pff}
This follows from Theorem~\ref{KFCQ2} and Corollary~\ref{Canc}.
\end{pff}

\section{Kakutani-Rokhlin decompositions}\label{Sec:KRd}

\indent
In this section, we show that the \tggp\  coming from
a free minimal action of $\Z^d$ on the Cantor set
is an almost AF Cantor groupoid.
The argument is essentially a reinterpretation of some of
the results of Forrest's paper~\cite{Fr} in terms of groupoids.
As can be seen by combining the results of this section and
Section~\ref{Sec:AAFG}, Forrest does not need all the conditions
he imposes to get his main theorem.
We use only the F{\o}lner condition, not the inradius conditions.

We presume that the construction of~\cite{Fr} can be generalized
to cover actions of more general groups.
Accordingly, we state the definitions in greater generality,
in particular using the word length metric rather than the
Euclidean distance on $\Z^d$.
This change has no effect on the results.

It should also be possible to generalize to actions that are
merely essentially free.
This generalization would require more substantial modification of
Forrest's definitions, and we do not carry it out here.
Doing so would gain no generality for the groups we actually handle,
namely $\Z^d$.
An essentially free minimal action of an abelian group is necessarily
free, because the fixed points for any one group element form a
closed invariant set.

We begin by establishing notation for this section.

\begin{cnv}\label{GpConv}
Throughout this section:

(1)
$\Gm$ is a countable discrete group
with a fixed finite generating set $\Sm$ which is symmetric in the
sense that $\gm \in \Sm$ implies $\gm^{-1} \in \Sm$.

(2)
$X$ is the Cantor set,
and $\Gm$ acts freely and \ct ly on $X$ on the left.
The action is denoted $(\gm, x) \mapsto \gm x$.

(3)
$G = \Gm \times X$ is the \tggp, and is equipped with
the Haar system consisting of counting measures.
Thus $G$ is a Cantor groupoid.
(See Definition~\ref{CGDfn} and Example~\ref{CTGEx}.)
\end{cnv}

For the convenience of the reader, we reproduce here the relevant
definitions from~\cite{Fr}, but stated for the more general situation
of Convention~\ref{GpConv}.
(Forrest considers the case in which $\Gm = \Z^d$ and $\Sm$ consists of
the standard basis vectors and their inverses.)

\begin{dfn}\label{KRd}
(Definition~2.1 of~\cite{Fr}.)
Let the notation be as in Convention~\ref{GpConv}.

(1)
A {\emph{tower}} is a pair $(E, S)$, in which $E \S X$ is a
compact open subset and $S \S \Gm$ is a finite subset, such that
$1 \in S$ and the sets $\gm E$, for $\gm \in S$, are pairwise disjoint.

(2)
The {\emph{levels}} of a tower $(E, S)$ are the sets
$\gm E \S X$ for $\gm \in S$.

(3)
A {\emph{traverse}} of a tower $(E, S)$ is a set of the
form $\{ \gm x \colon \gm \in S \}$ with $x \in E$.

(4)
A {\emph{\KRd}} $Q$ is a finite collection of towers whose
levels form a partition $P_Q$ of $X$
(called the {\emph{partition determined by $Q$}}).
\end{dfn}

In the following definition, we use the word length metric
(as opposed to the Euclidean norm on $\Z^d$ used in
Definition~3.2 of~\cite{Fr}).
However, in that case the two metrics are equivalent in the usual sense
for norms on Banach spaces, so the difference is not significant.


\begin{dfn}\label{Folner}
Let $l (\gm)$ (or $l_{\Sm} (\gm)$) denote the word length of
$\gm \in \Gm$ relative to the generating set $\Sm$.
Let $S \S \Gm$ be a nonempty finite subset.
The {\emph{F{\o}lner constant}} $c (S)$ (or $c_{\Sm} (S)$) is
the least number $c > 0$ such that
\[
\card (S \triangle \gm S) \leq c \cdot l (\gm) \card (S)
\]
for all $\gm \in \Gm$.
Here $S \triangle \gm S$ is the symmetric difference
\[
S \triangle \gm S = ( S \SM \gm S) \cup ( \gm S \SM S).
\]

%
The {\emph{F{\o}lner constant}} $c (Q)$ of a \KRd\  $Q$ is
\[
c (Q) = \sup_{(E, S) \in Q} c (S).
\]
\end{dfn}

\begin{dfn}\label{Refinement}
(Definition~2.2 of~\cite{Fr}.)
Let $Q_1$ and $Q_2$ be \KRd s.
We say that $Q_2$ {\emph{refines}} $Q_1$ if:
\begin{itemize}
\item[(1)]
The partition $P_{Q_2}$ (Definition~\ref{KRd}(4)) refines
the partition $P_{Q_1}$.
\item[(2)]
Every traverse (Definition~\ref{KRd}(3)) of a tower in $Q_2$
is a union of traverses of towers in $Q_1$.
\end{itemize}
\end{dfn}

\begin{dfn}\label{TSbGpd}
Let $(E, S)$ be a tower as in Definition~\ref{KRd}(1), and let
$G = \Gm \times X$ as in Convention~\ref{GpConv}(3).
We define the subset $G_{(E, S)} \S G$ by
\[
G_{(E, S)} = \{ ( \gm_1 \gm_2^{-1}, \, \gm_2 x )
   \colon x \in E \, \, {\mbox{and}} \, \, \gm_1, \gm_2 \in S \}.
\]
We show below that it is a subgroupoid of $G$.
We equip it with the Haar system consisting of counting measures;
we show below that this really is a Haar system.

Let $Q$ be a \KRd\  as in Definition~\ref{KRd}(4).
We define the {\emph{subgroupoid associated with $Q$}} to be
\[
G_Q = \bigcup_{(E, S) \in Q} G_{(E, S)} \S G.
\]
Again, we show below that it is a subgroupoid of $G$, and that we may
equip it with the Haar system consisting of counting measures.
\end{dfn}

\begin{lem}\label{SubgpdIsCptOpen}
Assume the hypotheses of Convention~\ref{GpConv}.
Then the subset
$G_{(E, S)} \S G$ of Definition~\ref{TSbGpd} is a compact open
Cantor subgroupoid of $G$.
Further, the subset $G_Q \S G$
of Definition~\ref{TSbGpd} is a compact open Cantor subgroupoid of $G$
which contains the unit space $G^{(0)}$ of $G$,
and is the disjoint union of the sets $G_{(E, S)}$ as $(E, S)$ runs
through the towers of $Q$.
\end{lem}

\begin{pff}
The subset $G_{(E, S)}$ is a finite union of subsets of $G$ of the
form
\[
\{ ( \gm_1 \gm_2^{-1}, \, \gm_2 x ) \colon x \in E \}
\]
for fixed $\gm_1, \gm_2 \in \Gm$.
Since $E$ is compact and open in $X$,
these sets are compact and open in $G$.
Therefore $G_{(E, S)}$ is compact and open in $G$.
So $G_Q$, being a finite union of sets of the form $G_{(E, S)}$,
is also compact and open in $G$.
It is immediate to check that the subset $G_{(E, S)}$ is closed under
product and inverse in the groupoid $G$, and so is a subgroupoid.
It is a Cantor groupoid by Example~\ref{OpenSbGpd}.

The subset $G_Q$ is the disjoint union of finitely many sets
of the form $G_{(E, S)}$, and is therefore also a Cantor groupoid.
Its unit space, thought of as a subset of $X$, is
\[
\bigcup_{(E, S) \in Q} G_{(E, S)}^{(0)}
 = \bigcup_{(E, S) \in Q} \left( \bigcup_{\gm \in S} \gm E \right),
\]
which is the union of all the levels of towers in $Q$, and is therefore
equal to $X$ (Definition~\ref{KRd}(4)).
Thus it is equal to $G^{(0)}$.
\end{pff}

\begin{lem}\label{GpdRef}
Assume the hypotheses of Convention~\ref{GpConv}.
Let $Q_1$ and $Q_2$ be \KRd s such that $Q_2$ refines $Q_1$
(Definition~\ref{Refinement}).
Then $G_{Q_1} \S G_{Q_2}$.
\end{lem}

\begin{pff}
Let $(E, S)$ be a tower of $Q_1$ and let $g \in G_{(E, S)}$.
Then $g = (\gm_1 \gm_2^{-1}, \, \gm_2 x)$ for some $x \in E$ and
$\gm_1, \, \gm_2 \in S$.
Since $P_{Q_2}$ is a partition of $X$,
there exist a tower $(F, T) \in Q_2$ and elements $\et_2 \in T$ and
$y \in F$ such that $\gm_2 x = \et_2 y$.

The traverse $\{ \et y \colon \et \in T \}$
(Definition~\ref{KRd}(3)) of $(F, T)$ is,
by Definition~\ref{Refinement}(2),
a union of traverses of towers in $Q_1$.
That is, there exist towers
\[
(E_1, S_1), \, \dots, \, (E_n, S_n) \in Q_1
\]
(not necessarily distinct)
and points $x_k \in E_k$ for $1 \leq k \leq n$, such that
\[
\{ \et y \colon \et \in T \}
   = \bigcup_{k = 1}^n \{ \gm x_k \colon \gm \in S_k \}.
\]
Therefore $\et_2 y = \gm x_k$ for some $k$ and some $\gm \in S_k$.
Now $\gm x_k$ is in the level $\gm E_k \in P_{Q_1}$.
Since $\gm x_k = \gm_2 x$, the action is free, and the levels of
all the towers of $Q_1$ form a partition of $X$, it follows that
\[
(E_k, S_k) = (E, S), \,\,\,\,\,\, x_k = x, \andeqn \gm = \gm_2.
\]
In particular,
\[
\{ \gm x \colon \gm \in S \} \S \{ \et y \colon \et \in T \},
\]
so there is $\et_1 \in T$ such that $\gm_1 x = \et_1 y$.
Substituting $y = \et_2^{-1} \gm_2 x$, we get
$\gm_1 x = \et_1 \et_2^{-1} \gm_2 x$.
Because the action is free,
it follows that $\et_1 \et_2^{-1} = \gm_1 \gm_2^{-1}$.
We now have
\[
g = (\gm_1 \gm_2^{-1}, \, \gm_2 x)
  = (\et_1 \et_2^{-1}, \, \et_2 y) \in G_{(F, T)} \S G_{Q_2}.
\]
\end{pff}

\begin{thm}\label{ExistRef}
(Forrest~\cite{Fr}.)
In the situation of Convention~\ref{GpConv}, assume that $\Gm = \Z^d$,
that $\Sm$ consists of the standard basis vectors and their inverses,
and that the action of $\Z^d$ on $X$ is minimal (as well as being free).
Then there exists a sequence
\[
Q_1, \, Q_2, \, Q_3, \, \dots
\]
of \KRd s such that $Q_{n + 1}$ refines $Q_n$ for all $n$, and such that
the F{\o}lner constants obey $\limi{n} c (Q_n) = 0$.
\end{thm}

\begin{pff}
Combine Proposition~5.1 of~\cite{Fr} with Lemma~3.1 of~\cite{Fr}.
This immediately gives everything except for the estimate
on the F{\o}lner constants.
{}From~\cite{Fr} we get an estimate on constants slightly different
from those given in Definition~\ref{Folner}.
Specifically, define $c_0 (S)$ and $c_0 (Q)$ as in
Definition~\ref{Folner}, but considering $\Z^d$ as a subset of $\R^d$
and using $\| \gm \|_2$ in place of $l ( \gm )$.
That is,
\[
\card (S \triangle \gm S) \leq c_0 (S) \cdot \| \gm \|_2 \card (S),
\]
etc.
That $\limi{n} c_0 (Q_n) = 0$ follows from the fact that the atoms of
the partitions ${\mathcal{T}}_n (x)$ in Proposition~5.1 of~\cite{Fr},
for fixed $n$ and as $x$ runs through $X$,
are exactly the translates in $\Gm$
of the sets $S$ appearing in the towers $(E, S)$ which make up $Q_n$.
Now for any $\gm \in \Z^d$ we have
$l (\gm) = \| \gm \|_{\infty} \geq d^{-1} \| \gm \|_2$.
Therefore $c (Q_n) \leq d \cdot c_0 (Q_n)$ for all $n$,
whence $\limi{n} c (Q_n) = 0$, as desired.
\end{pff}

\begin{thm}\label{TGIsAF}
Assume the hypotheses of Convention~\ref{GpConv},
and assume moreover that
the action of $\Gm$ on $X$ is minimal (as well as being free).
Assume that there exists a sequence
\[
Q_1, \, Q_2, \, Q_3, \, \dots
\]
of \KRd s such that $Q_{n + 1}$ refines $Q_n$ for all $n$, and such that
the F{\o}lner constants obey $\limi{n} c (Q_n) = 0$.
Then $G = \Gm \times X$ is almost AF in the sense of
Definition~\ref{TAFGDfn}, and $C_{\mathrm{r}}^* (G)$ is simple.
\end{thm}

\begin{pff}
The algebra $C_{\mathrm{r}}^* (G)$ is the same as the reduced
\tgca\  $C_{\mathrm{r}}^* (\Gm, X)$ by Proposition~\ref{CrPisGpC}.
So it is simple by the corollary at the end of~\cite{AS}.
(See the preceding discussion and Definition~1 there.)
To prove that $G$ is almost AF, it is by
Proposition~\ref{SimpleAAF} now enough to verify
Definition~\ref{TAFGDfn}(1).

With $G_{Q_n}$ as in Definition~\ref{TSbGpd},
define $G_0 = \bigcup_{n = 1}^{\I} G_{Q_n}$.
Using Lemma~\ref{SubgpdIsCptOpen} and Lemma~\ref{GpdRef},
we easily verify the conditions of Definition~\ref{AFGDfn}.
Thus $G_0$ is an AF Cantor groupoid, which is open in $G$
because each $G_{Q_n}$ is open in $G$.

It remains to verify the condition that $s (K)$ be thin in $G_0$
(Definition~\ref{Thin}) for every compact set $K \S G \SM G_0$.
For this argument,
we will identify the unit space $G^{(0)} = \{ 1 \} \times X$ with $X$
in the obvious way.
So let $K \S G \SM G_0$ be compact and let $n \in \N$.
Let
\[
T = \{ \gm \in \Gm
      \colon (\{ \gm \} \times X) \cap K \neq \varnothing \}
  \S \Gm.
\]
Then $T$ is finite because $K$ is compact.
Let $l$ be the word length
metric on $\Gm$, as used in Definition~\ref{Folner}.
Let $\rh = \sup_{\gm \in T} l (\gm)$.
Choose $m$ so large that
\[
c (Q_m) < \frac{1}{\rh n \cdot \card (T)}.
\]

We now claim that if $(E, S) \in Q_m$ is a tower and
\[
x \in s (K) \cap \bigcup_{\et \in S} \et E,
\]
then there exist $\gm \in T$ and $\et \in S \SM \gm^{-1} S$ such that
$(\gm, x) \in \{ \gm \} \times \et E$.
To prove this, recall that $s (\gm, x) = x$ (really $(1, x)$),
and find $\gm \in T$ such that $(\gm, x) \in K$.
Also choose $\et \in S$ such that $\et^{-1} x \in E$.
Then write $(\gm, x) = ( [\gm \et] \et^{-1}, \, \et [ \et^{-1} x ] )$.
Since $(\gm, x) \not\in G_{(E, S)}$, we have $\gm \et \not\in S$.
So $\et \in S \SM \gm^{-1} S$.
This proves the claim.

It follows that the sets $\et E$, for $(E, S) \in Q_m$, $\gm \in T$, and
$\et \in S \SM \gm^{-1} S$, form a cover of $s (K)$ by finitely many
disjoint compact open subsets of $X$.
For
\[
(E, S) \in Q_m \andeqn \et \in \bigcup_{\gm \in T} (S \SM \gm^{-1} S),
\]
define $L_{(E, S), \et} = s (K) \cap \et E$.
The sets $L_{(E, S), \et}$ are a cover of $s (K)$ by finitely many
disjoint compact subsets of $X$.

For fixed $(E, S) \in Q_m$, we now claim that there
are compact $G_{(E, S)}$-sets
\[
M_{(E, S), 1}, \, M_{(E, S), 2},
  \, \dots, \, M_{(E, S), n} \S G_{(E, S)}
\]
such that
\[
s (M_{(E, S), k}) = s (K) \cap \bigcup_{\et \in S} \et E
\]
and the sets
$r (M_{(E, S), 1}), \, r (M_{(E, S), 2}),
  \, \dots, \, r (M_{(E, S), n})$
are pairwise disjoint.
Since the subgroupoids $G_{(E, S)}$
are disjoint (Lemma~\ref{SubgpdIsCptOpen}) and their union
is contained in $G_0$, we can take then the required $G_0$-sets to be
\[
M_k = \bigcup_{(E, S) \in Q_m} M_{(E, S), k}.
\]

To find the $M_{(E, S), k}$, first note that for $\gm \in T$, we have
$l (\gm^{-1}) = l (\gm) \leq \rh$, so, using Definition~\ref{Folner},
\[
\card (S \SM \gm^{-1} S) \leq \card (S \triangle \gm^{-1} S)
 \leq c (Q_m) \rh \cdot \card (S) < \frac{\card (S)}{n \cdot \card (T)}.
\]
Therefore, with
\[
R_{(E, S)} = \bigcup_{\gm \in T} (S \SM \gm^{-1} S),
\]
we have
$\card ( R_{(E, S)} ) < \frac{1}{n} \card (S)$.
So there exist $n$ injective functions
\[
\sm_1, \, \sm_2, \, \dots, \, \sm_n : R_{(E, S)} \to S
\]
with disjoint ranges.
Define
\[
M_{(E, S), k} = \bigcup_{\et \in R_{(E, S)} }
 \left( \{ \sm_k (\et) \et^{-1} \} \times L_{(E, S), \et} \right).
\]
Then $M_{(E, S), k} \S G_{(E, S)}$
because $L_{(E, S), \et} \S \et E$.
Clearly $M_{(E, S), k}$ is compact.
The restriction of the source map to $M_{(E, S), k}$ is injective
because the $L_{(E, S), \et}$ are disjoint.
Further, $\sm_k (\et) \et^{-1} L_{(E, S), \et} \S \sm_k (\et) E$,
so the sets $\sm_k (\et) \et^{-1} L_{(E, S), \et}$,
for $\et \in R_{(E, S)}$ and $1 \leq k \leq n$, are pairwise disjoint.
It follows both that the restriction of the range map to
$M_{(E, S), k}$ is injective and that the sets
\[
r (M_{(E, S), 1}), \, r (M_{(E, S), 2}),
                    \, \dots, \, r (M_{(E, S), n})
\]
are pairwise disjoint.
Thus the $M_{(E, S), k}$
are $G_{(E, S)}$-sets with the required properties.
We have verified Definition~\ref{Thin} for $s (K)$.
\end{pff}

\begin{cor}\label{ZdisAAF}
Let $d$ be a positive integer, let $X$ be the Cantor set, and let
$\Z^d$ act freely and minimally on $X$.
Then the \tggp\  $G = \Z^d \times X$ is almost AF in the sense of
Definition~\ref{TAFGDfn}, and $C_{\mathrm{r}}^* (G)$ is simple.
\end{cor}

\begin{pff}
This is immediate from Theorems~\ref{ExistRef} and~\ref{TGIsAF}.
\end{pff}

\begin{thm}\label{MainZd}
Let $d$ be a positive integer, let $X$ be the Cantor set, and let
$\Z^d$ act freely and minimally on $X$.
Then:
\begin{itemize}
\item[(1)]
$C^* (\Z^d, X)$ has (topological) stable rank one
in the sense of~\cite{Rf}.
\item[(2)]
$C^* (\Z^d, X)$ has real rank zero in the sense of~\cite{BP}.
\item[(3)]
Let $p, q \in \Mi (C^* (\Z^d, X))$
be \pj s such that $\ta (p) < \ta (q)$
for all normalized traces $\ta$ on $C^* (\Z^d, X)$.
Then $p$ is \mvnt\  to a sub\pj\  of $q$.
\end{itemize}
\end{thm}

\begin{pff}
Using Proposition~\ref{CrPisGpC} and Corollary~\ref{ZdisAAF},
we obtain Part~(1) from Theorem~\ref{tsr1},
Part~(2) from Theorem~\ref{RR0}, and
Part~(3) from Corollary~\ref{FCQ2}.
\end{pff}

\section{Tilings and quasicrystals}\label{Sec:QC}

\indent
In this section, we show that the \ca s associated with three different
broad classes of aperiodic tilings have real rank zero and
stable rank one, and satisfy Blackadar's Second
Fundamental Comparability Question.
In particular,
we strengthen the conclusion of the main result of~\cite{Pt5}.
Then we discuss the relationship with the Bethe-Sommerfeld Conjecture
for quasicrystals.

We begin by showing that the groupoids of~\cite{Pt5} are almost AF.
The proof consists of assembling, in the right order, various
results proved in~\cite{Pt5}.

\begin{thm}\label{PtIsAAF}
Consider a substitution tiling system in $\R^d$ as in Section~1
of~\cite{Pt5}, under the conditions imposed there:
\begin{itemize}
\item
The substitution is primitive.
\item
The finite pattern condition holds.
\item
The system is aperiodic.
\item
The substitution forces its border.
\item
The capacity of the boundary of every prototile is strictly
less than $d$.
\end{itemize}
Then the groupoid $R_{\text{punc}}$ defined in Section~1
of~\cite{Pt5} is an almost AF Cantor groupoid.
\end{thm}

\begin{pff}
That the groupoid $R_{\text{punc}}$ is a Cantor groupoid in the
sense of Definition~\ref{CGDfn} follows from the construction as
described on pages 594--595 of~\cite{Pt5} and on page~187 of~\cite{KP}.
Note that the base for the topology described in~\cite{Pt5} is
countable, so that $R_{\text{punc}}$ really is second countable.

We take the open AF subgroupoid required in Definition~\ref{TAFGDfn}
to be the
open subgroupoid $R_{\text{AF}}$ from page~596 of~\cite{Pt5}.
It is easily seen from the construction there to be an AF Cantor
groupoid which contains the unit space of $R_{\text{punc}}$.
(Also see pages 198--200 of~\cite{KP}.)

We must show that if $K \subset R_{\text{punc}} \setminus R_{\text{AF}}$
is compact, then $s (K)$ is thin in the sense of Definition~\ref{Thin}.
Referring to the definition of $R_{\text{punc}}$ and its topology
(pages 594--595 of~\cite{Pt5}), we see that $R_{\text{punc}}$
consists of certain pairs $(T, \, T + x)$ in which $T$ is a tiling
of $\R^d$ and $x \in \R^d$, and that there is a base for its topology
consisting of sets of the form $\{ (T, \, T + x) \colon T \in U \}$
for suitable $x \in \R^d$ and suitable sets $U$ of tilings.
It is immediate from this that for any $r > 0$ the set
\[
\{ (T, \, T - x) \in R_{\text{punc}} \colon \| x \| < r \}
\]
is open in $R_{\text{punc}}$.
As $r$ varies, these sets cover $R_{\text{punc}}$, and $K^{-1}$ is
also a compact subset of $R_{\text{punc}} \setminus R_{\text{AF}}$,
so there is $r > 0$ such that $K^{-1}$ is contained in the set
\[
L = \{ (T, \, T - x) \in R_{\text{punc}} \setminus R_{\text{AF}} \colon
           \| x \| \leq r \}.
\]
Since $s (K) \subset r (L)$, it suffices to show that $r (L)$ is thin.

The proof of Theorem~2.1 of~\cite{Pt5}, at the end of Section~2.1 there,
consists of showing that there are compact open sets
$U_n \subset \Om_{\text{punc}}$, the unit space of $R_{\text{punc}}$,
for $n \in \N$, each containing $r (L)$, and \hme s $\gm_n$ from
$U_n$ to disjoint subsets of $\Om_{\text{punc}}$, each having graph
contained in $R_{\text{AF}}$.
But it is immediate from this that $r (L)$ is thin.
This completes the verification of Part~(1) of Definition~\ref{TAFGDfn}.

{}From the description of $C^* ( R_{\text{AF}} )$ and its K-theory at
the end of Section~1 of~\cite{Pt5}, we see that this algebra is a direct
limit of a system of finite dimensional \ca s in which the matrix
$B$ of partial embedding multiplicities is the same at each stage.
Moreover, there is $n$ such that $B^n$ has no zero entries.
Therefore the direct limit algebra is simple, whence
$C^*_{\mathrm{r}} ( R_{\text{AF}} )$ is simple.
It now follows from Proposition~\ref{SimpleAAF}
that $R_{\text{punc}}$ is almost AF.
\end{pff}

\smallskip

The following corollary contains Theorem~1.1 of~\cite{Pt5}.
Using more recent results (see Theorem~5.1 of \cite{KmP}),
it is now also possible to obtain this result from Theorem~\ref{MainZd}
in the same way the that next two theorems are proved.

\begin{cor}\label{PtCSAlg}
For a substitution tiling system in $\R^d$ as in Theorem~\ref{PtIsAAF},
the \ca s of the associated groupoid $R_{\text{punc}}$
as in~\cite{Pt5} has stable
rank one and real rank zero, and satisfies Blackadar's Second
Fundamental Comparability Question (\cite{Bl4}, 1.3.1).
\end{cor}

\begin{pff}
Using Theorem~\ref{PtIsAAF},
we obtain stable rank one from Theorem~\ref{tsr1},
real rank zero from Theorem~\ref{RR0}, and
Blackadar's Second Fundamental
Comparability Question from Corollary~\ref{FCQ2}.
\end{pff}

\smallskip

We can also obtain the same result for several other kinds of
aperiodic tilings.
As discussed on page~198 of~\cite{KP},
we reduce to the case of crossed products
by free minimal actions of $\Z^d$ on the Cantor set.

\begin{thm}\label{PjMthTiling}
Consider a projection method pattern ${\mathcal{T}}$ as in~\cite{FHK},
with data $(E, K, u)$ (see Definitions I.2.1 and I.4.4 of~\cite{FHK})
such that $E \cap \Z^N = \{ 0 \}$.
Let ${\mathcal{G}}{\mathcal{T}}$ be the associated groupoid
(Definition II.2.7 of~\cite{FHK}).
Then $C^* ( {\mathcal{G}}{\mathcal{T}} )$ has stable
rank one and real rank zero, and satisfies Blackadar's Second
Fundamental Comparability Question.
\end{thm}

\begin{pff}
It follows from Theorem II.2.9 and Corollary I.10.10 of~\cite{FHK} that
$C^* ( {\mathcal{G}}{\mathcal{T}} )$ is strongly Morita equivalent to
the \tgca\  of a minimal action of $\Z^d$ on a Cantor set
$X_{\mathcal{T}}$, for a suitable $d$.
The action is free since the action of $\R^d$ in the pattern dynamical
system in that corollary is free by construction.

The required properties hold for $C^* (\Z^d, \, X_{\mathcal{T}})$
by Theorem~\ref{MainZd}.
Since $C^* ( {\mathcal{G}}{\mathcal{T}} )$ is strongly Morita equivalent
to $C^* (\Z^d, \, X_{\mathcal{T}})$, and both algebras are separable,
these two algebras are stably isomorphic by Theorem~1.2 of~\cite{BGR}.
So $C^* ( {\mathcal{G}}{\mathcal{T}} )$
has stable rank one by Theorem~3.6 of~\cite{Rf} and has
real rank zero by Theorem~3.8 of~\cite{BP}.
It is clear that Blackadar's Second Fundamental
Comparability Question is preserved by stable isomorphism.
\end{pff}

\smallskip

Theorem II.2.9 of~\cite{FHK} shows that several other groupoids
associated to ${\mathcal{T}}$ give \ca s strongly Morita equivalent
to $C^* ( {\mathcal{G}}{\mathcal{T}} )$.
These \ca s therefore also have stable
rank one and real rank zero, and satisfy Blackadar's Second
Fundamental Comparability Question.

\begin{thm}\label{GenGridT}
Consider a minimal aperiodic
generalized dual (or grid) method tiling~\cite{SSL},
satisfying the condition~(D2) before Lemma~10 of~\cite{Kl}.
Then the \ca\  associated with the tiling has stable
rank one and real rank zero, and satisfies Blackadar's Second
Fundamental Comparability Question.
\end{thm}

\begin{pff}
Theorem~1 and Lemma~11 of~\cite{Kl} show that the \ca\  associated
with such a tiling is stably isomorphic to the crossed product
\ca\  for an action of $\Z^d$ on the Cantor set.
(See Definition~12 of~\cite{Kl} and the following discussion,
the discussion following Theorem~1 of~\cite{Kl},
and also Definition~13 of~\cite{Kl}.)
The action is minimal
(see the discussion following Definition~12 of~\cite{Kl})
and free (see Definition~8 and Lemma~3 of~\cite{Kl}).
Therefore the crossed product has stable
rank one and real rank zero, and satisfies Blackadar's Second
Fundamental Comparability Question, by Theorem~\ref{MainZd}.
It follows as in the proof of Theorem~\ref{PjMthTiling} that
the \ca\  of the tiling has these properties as well.
\end{pff}

\smallskip

We now turn to the Bethe-Sommerfeld conjecture for quasicrystals.
Consider the Schr\"{o}dinger operator $H$ for an electron moving in a
solid in $\R^d$, which at this point may be crystalline, amorphous, or
quasicrystalline.
(The physical case is, of course, $d = 3$.)
Associated with this situation there is a \ca\  which contains
all bounded \cfn s of $H$ and all its translates by elements of $\R^d$.
See~\cite{BHZ}.
For crystals, in which the locations of the atomic nuclei form
a lattice in $\R^d$, the Bethe-Sommerfeld conjecture asserts that
if $d \geq 2$ then
there is an energy above which the spectrum of $H$ has no gaps.
See~\cite{DT},~\cite{Skr}, Corollary~2.3 of~\cite{HM}, and~\cite{Krp}
for some of the mathematical results on this conjecture.
According to the introduction to~\cite{BHZ}, it is expected that
this property also holds for quasicrystals.
We point out that quasicrystals actually occur in nature;
see~\cite{SG} for a survey.
For the Schr\"{o}dinger operator in the so-called tight binding
representation for a quasicrystal whose atomic sites are given by
an aperiodic tiling, the relevant \ca\  is the \ca\  of the tiling.
See Section~4 of~\cite{BHZ}.

We have just seen that the \ca s associated with three different broad
classes of tilings have real rank zero, that is, every selfadjoint
element, including $H$, can be approximated arbitrarily closely in
norm by selfadjoint elements with finite spectrum.
In Proposition~\ref{TDSpIsGen} below, we further show that among the
selfadjoint elements of such an algebra, those with totally disconnected
spectrum are generic, in the sense that they form a dense $G_{\dt}$-set.
That is, the selfadjoint elements of the relevant \ca\  whose
spectrum is not totally disconnected form a meager set in the
sense of Baire category.

These results say nothing about the spectrum of $H$ itself.
It is also not clear whether arbitrary norm small selfadjoint
perturbations are physically relevant.
Nevertheless, gap labelling theory (see~\cite{BHZ} for a recent survey)
exists in arbitrary dimensions,
and the introduction of~\cite{BHZ} raises the question of its
physical interpretation when there are no gaps.
Our results do suggest the possibility that this theory has
physical significance in the general situation.

\begin{prp}\label{TDSpIsGen}
Let $A$ be a \ca\  with real rank zero.
Then there is a dense $G_{\delta}$-set $S$ in the selfadjoint part
$A_{\text{sa}}$ of $A$ such that every element of $S$ has
totally disconnected spectrum.
\end{prp}

\begin{pff}
For a finite subset $F \subset \R$, let
\[
S_F = \{ a \in A_{\text{sa}} \colon \spec (a) \cap F = \varnothing \}.
\]
Then $S_F$ is open $A_{\text{sa}}$,
since if $a, b \in A_{\text{sa}}$ satisfy $\|a - b\| < \ep$,
then $\spec (b)$
is contained in the $\ep$-neighborhood of $\spec (a)$.
We show that $S_F$ is dense in $A_{\text{sa}}$.
Given any selfadjoint element $a \in A$ and $\ep > 0$, the real
rank zero condition provides a selfadjoint element $b \in A$ with
finite spectrum and $\|a - b\| < \frac{1}{2} \ep$.
By perturbing eigenvalues, it is easy to find a selfadjoint element
$c \in A$ with $\|b - c\| <  \frac{1}{2} \ep$ and such that
$\spec (c) \cap F = \varnothing$.
Then $c \in S_F$ and $\|a - c\| < \ep$.
So $S_F$ is dense.

The set of selfadjoint elements $a \in A$ with
$\spec (a) \cap \Q = \varnothing$ is a countable intersection of
sets of the form $S_F$, and is therefore a dense $G_{\delta}$-set.
Clearly all its elements have totally disconnected spectrum.
\end{pff}

\section{Examples and open problems}\label{Sec:Problems}

\indent
In this section, we discuss some open problems, beginning with
the various possibilities for improving our results.
One of the most obvious questions is the following:

\begin{qst}
Let $G$ be an almost AF Cantor groupoid.
Does it follow that $C_{\mathrm{r}}^* (G)$
is tracially AF in the sense of Definition~2.1 of~\cite{Ln10}?
\end{qst}

Lemma~\ref{CutToAF} and the thinness condition in
Definition~\ref{TAFGDfn} come close to giving the tracially AF
property;
the main part that is missing is the requirement that the
projections one chooses approximately commute with the
elements of the finite set.

By Theorem~3.1 of~\cite{FH}, the K-theory of the crossed
product of the Cantor set by any action of $\Z^d$ is torsion free.
Therefore, by recent work of H.\  Lin~\cite{Ln15}, a positive
answer to this question
would imply a positive solution to the following conjecture:

\begin{cnj}
Let $d$ be a positive integer, let $X$ be the Cantor set, and let
$\Z^d$ act freely and minimally on $X$.
Then $C^* (\Z^d, X)$ is an AT algebra, that is,
isomorphic to a direct limit $\dirlim A_n$ of \ca s $A_n$ each of
which is a finite direct sum of \ca s of the form $M_k$ or
$C (S^1, \, M_k)$.
(The matrix sizes $k$ may vary among the summands, even
for the same value of $n$.)
\end{cnj}

This conjecture is also predicted by the Elliott classification
conjecture.
It is true for $d = 1$ (essentially~\cite{Pt2}).
We point out here that Matui has shown in~\cite{Mt}
(see Proposition~4.1 and Theorem~4.8) that when $d = 2$, the
algebra $C^* (\Z^d, X)$ is often AF embeddable.

We next address the possibility of generalizing the group.
For not necessarily abelian groups, we should presumably
require that the action be minimal and merely essentially free,
which for minimal actions means that the set of $x \in X$ with trivial
isotropy subgroup is dense in $X$.
This is the condition needed to ensure simplicity of the \tgca.
See the final corollary of~\cite{AS}, noting that all actions of
amenable groups are regular (as described before this corollary)
and that essential freeness of an action of a countable discrete
group on a compact metric space $X$ is equivalent to
topological freeness, Definition~1 of~\cite{AS}, of the
action on $C (X)$ (as is easy to prove).
Thus, let $\Gm$ act minimally and essentially freely on the
Cantor set $X$.
If $\Gm$ is close to $\Z^d$, and especially if the action is
actually free, then there seems to be
a good chance of adapting the methods of Forrest~\cite{Fr}
to show that the \tggp\  is still almost AF.
This might work, for example, if $\Gm$ has a finite index subgroup
isomorphic to $\Z^d$.
In this context, we mention that crossed products of
free minimal actions of the free product
$\Z / 2 \Z \bigstar \Z / 2 \Z$ on the
Cantor set, satisfying an additional technical condition,
are known to be AT, in fact AF~\cite{BEK}.
The group $\Z / 2 \Z \bigstar \Z / 2 \Z$ has an index two
subgroup isomorphic to $\Z$.
However, we believe the correct generality is to allow $\Gm$ to
be an arbitrary countable amenable group.

\begin{qst}
Let $X$ be the Cantor set, and let the countable amenable group
$\Gm$ act minimally and essentially freely on $X$.
Is the \tggp\  $\Gm \times X$ almost AF?
\end{qst}

Even if not, do the conclusions of Theorem~\ref{MainZd} still hold?
That is:

\begin{qst}
Let $X$ be the Cantor set, and let the countable amenable group
$\Gm$ act minimally and essentially freely on $X$.
Does it follow that $C^* (\Gm, X)$ has stable rank one and
real rank zero?
If $p, q \in \Mi (C^* (\Z^d, X))$
are \pj s such that $\ta (p) < \ta (q)$
for all normalized traces $\ta$ on $C^* (\Z^d, X)$, does it follow that
$p$ is \mvnt\  to a sub\pj\  of $q$?
\end{qst}

The other obvious change is to relax the condition on the space.
If $X$ is not totally disconnected, then even a crossed product by
a free minimal action of $\Z$ need not have real rank zero.
Examples~4 and~5 in Section~5 of~\cite{Cn} have no nontrivial \pj s.
However, it seems reasonable to hope that the other parts of the
conclusion of Theorem~\ref{MainZd} still hold, at least under mild
restrictions.

\begin{qst}
Let $X$ be a compact metric space with finite covering dimension
\cite{En},
and let the countable amenable group
$\Gm$ act minimally and essentially freely on $X$.
Does it follow that $C^* (\Gm, X)$ has stable rank one?
If $p, q \in \Mi (C^* (\Z^d, X))$
are \pj s such that $\ta (p) < \ta (q)$
for all normalized traces $\ta$ on $C^* (\Z^d, X)$, does it follow that
$p$ is \mvnt\  to a sub\pj\  of $q$?
\end{qst}

We think we know how to prove stable rank one when $\Gm = \Z^d$,
but a number of technical details need to be worked out.

The definition given for an almost AF Cantor groupoid whose reduced
\ca\  is not simple is merely a guess, and it is also not clear
what properties the reduced \ca\  of such a groupoid should have.

\begin{qst}
What is the ``right'' definition of an almost AF Cantor groupoid
in the case that the reduced \ca\  is not simple?
\end{qst}

\begin{qst}
Let $G$ be an almost AF Cantor groupoid such that
$C_{\mathrm{r}}^* (G)$ is not simple.
What structural consequences does this have for $C_{\mathrm{r}}^* (G)$?
\end{qst}

The definition should surely exclude the \tggp\  coming from the
action of $\Z$ on $\Z \cup \{ \pm \infty \}$ by translation,
even though this action is free on a dense set.
It is not essentially free, which is equivalent to the
\tggp\  $G$ not being essentially principal (Definition~\ref{PDfn}(4)).

Suppose, however, that $h_1, \, h_2 \colon X \to X$
are two commuting \hme s of the Cantor set $X$, such that the map
$(n_1, n_2) \mapsto h_1^{n_1} \circ h_2^{n_2}$ determines a free
minimal action of $\Z^2$.
Then the \tggp\  for the action of $\Z$ on $X$ generated by $h_1$
should be almost AF.
It is possible for this action to be already minimal,
or for $X$ to be the disjoint union of finitely many minimal sets,
in which case the situation is clear.
It is possible that neither of these happens---consider the product
of two minimal actions.
It is not clear what additional structure the \hme\  $h_1$ must
have (although it can't be completely arbitrary---see below),
and it is also not clear whether the crossed product must necessarily
have real rank zero, stable rank one, or order on its $K_0$-group
determined by traces.

Let $h$ be an arbitrary aperiodic \hme\  of the Cantor set $X$.
(That is, the action of $\Z$ it generates is free, but need not be
minimal.)
Theorem~3.1 of~\cite{Pt2} shows that if $h$ has no nontrivial
invariant subsets which are both closed and open, and if $h$ has
more than one minimal set, then $C^* (\Z, X, h)$ does not have
stable rank one.
If $h_1$ and $h_2$ are as above, and $h = h_1$, then the existence
of a unique minimal set $K$ for $h$ implies that $h$ is minimal.
(The set $K$ is necessarily invariant under $h_2$.)
On the other hand, we know of no examples of $h_1$ and $h_2$ as above
in which $h_1$ is not minimal yet has no nontrivial
invariant subsets which are both closed and open.

It follows from~\cite{Ps} that if
$h$ is an arbitrary aperiodic \hme\  of the Cantor set $X$,
then $C^* (\Z, X, h)$ has the ideal property,
that is, every ideal is generated as an ideal in the algebra by
its \pj s.
This certainly suggests that the reduced C*-algebra of an almost
AF Cantor groupoid should have the ideal property.
This condition is, however, rather weak in this context,
since every simple unital \ca, regardless of its real or stable
rank, has the ideal property.

We now give an example of an aperiodic \hme\  of the Cantor set
whose \tggp\  is not almost AF, but only because of the
failure of condition~(2) in Definition~\ref{TAFGDfn}.
Its \ca\  does not have stable rank one.

\begin{exa}\label{NonMinEx}
Let $X_1 = \Z \cup \{ \pm \infty \}$ be the two point compactification
of $\Z$, and let $h_1 \colon X_1 \to X_1$ be translation to the right
(fixing $\pm \infty$).
Let $X_2$ be the Cantor set, and let $h_2 \colon X_2 \to X_2$ be any
minimal \hme.
Let $X = X_1 \times X_2$, and let $h = h_1 \times h_2 \colon X \to X$.
Let $G = \Z \times X$ be the \tggp, as in Example~\ref{CTGEx},
for the corresponding action of $\Z$.
Then the $G$ has the following properties:
\begin{itemize}
\item
It is a Cantor groupoid.
\item
It satisfies condition~(1) of Definition~\ref{TAFGDfn}, the main part
of the definition of an almost AF Cantor groupoid.
\item
There is a nonempty open subset $U \subset G^{(0)}$ which is null for
all $G$-invariant Borel probability measures.
\item
The reduced \ca\  $C_{\mathrm{r}}^* (G)$ does not have stable rank one.
\end{itemize}

For the first statement,
observe that $X$ is a totally disconnected compact metric space with no
isolated points, and hence homeomorphic to the Cantor set.
So $G$ is a Cantor groupoid, as in Example~\ref{CTGEx}.

The third statement is also easy to prove, with
$U = \Z \times X_2$.
The $G$-invariant Borel probability measures on $G^{(0)}$ are exactly
the $h$-invariant Borel probability measures on $X$.
(We have not found this statement in the literature.
However, our case follows immediately from Remark~\ref{CGRmk}(2) and
the formulas in Example~\ref{CTGEx} by considering functions supported
on sets of the form $\{ \gm \} \times X$ and letting $\gm$ run through
the group.)
So let $\mu$ be any $h$-invariant Borel probability measure on $X$.
Then $\mu ( \{ n \} \times X_2)$ is independent of $n$,
and
\[
\sum_{n \in \Z} \mu ( \{ n \} \times X_2) = \mu (\Z \times X_2)
   \leq \mu (X) = 1.
\]
It follows that $\mu ( \{ n \} \times X_2) = 0$ for all $n \in \Z$,
whence $\mu (\Z \times X_2) = 0$.

We verify condition~(1) of Definition~\ref{TAFGDfn}.
Fix any $z \in X_2$.
Set $Y = X_1 \times \{ z \} \subset X$.
Let $G_0$ be the subgroupoid $G'$ as in Example~2.6 of~\cite{Pt4}
corresponding to this situation.
(Also see the statement of Theorem~2.4 of~\cite{Pt4}.)
Note that, in our notation, $G_0 = G \setminus (L \cup L^{-1})$, with
\[
L = \{ (k, \, h^l (y)) \colon {\mbox{$y \in Y$, $k, l \in \Z$,
     $l > 0$, and $k + l \leq 0$}} \}
   \subset G.
\]
To see that $G_0$ is AF (compare with~\cite{Pt1}), we choose a
decreasing sequence of compact open subsets $Z_n \subset X_2$ with
$\bigcap_{n = 1}^{\infty} Z_n = \{ z \}$.
Set $Y_n = X_1 \times Z_n$, define $L_n \subset G$ by using
$Y_n$ in place of $Y$ in the definition of $L$, and set
$H_n = G \setminus (L_n \cup L_n^{-1})$.
It is easy to check that $H_n$ is a closed and open subgroupoid of $G$
which is contained in $G_0$,
and that $G_0 = \bigcup_{n = 1}^{\infty} H_n$.
Since $h_2$ is minimal, there is $l (n)$ such that the sets
$h_2 (Z_n), \, h_2^2 (Z_n), \, \dots, \, h_2^{n (l)} (Z_n)$ cover $X_2$,
and it follows that
\[
H_n \subset \{ -n (l), \, - n (l) + 1, \, \dots, \, n (l) \} \times X.
\]
So $H_n$ is compact.
We have now shown that $G_0$ is in fact an AF Cantor groupoid.

We now show that if $K \subset G \setminus G_0$ is compact,
then $s (K)$ is thin in $G_0$.
We will identify the unit space $G^{(0)} = \{ 1 \} \times X$ with $X$
in the obvious way.
Compact subsets of thin sets are easily seen to be thin,
so it suffices to consider sets of the form
$K = (G \setminus G_0) \cap T$ with
\[
T = \{ - m, \, - m + 1, \, \dots, \, m - 1, \, m \} \times X
\]
for some $m \in \N$.
Now
\[
(G \setminus G_0) \cap T = ( L \cup L^{-1}) \cap T,
\]
so, following Example~\ref{CTGEx},
\[
s (K) = s (L \cap T) \cup r (L \cap T).
\]
Next, we calculate
\begin{align*}
L \cap T & =
  \{ (k, \, h^l (y)) \colon {\mbox{$y \in Y$, $k, l \in \Z$,
     $l > 0$, $k + l \leq 0$, and $-m \leq k \leq m$}} \}   \\
    & =
  \{ (- k, \, h^l (y)) \colon {\mbox{$y \in Y$, $k, l \in \Z$,
     and $1 \leq l \leq k \leq m$}} \}.
\end{align*}
It follows that
\[
s (L \cap T) = h (Y) \cup h^2 (Y) \cup \cdots \cup h^m (Y)
\]
and
\[
r (L \cap T) = \bigcup_{1 \leq l \leq k \leq m} h^{- k + l} (Y)
  = Y \cup h^{-1} (Y) \cup \cdots \cup h^{-m + 1} (Y).
\]

Now let $n \in \N$.
We find compact $G$-sets $S_1, S_2, \dots, S_n \S G$
such that $s (S_k) = s (K)$ and the sets
$r (S_1), \, r ( S_2), \, \dots, \, r (S_n)$ are pairwise disjoint.
For $1 \leq j \leq n$, define
\begin{align*}
S_j & = \{ m (j - 1) \}
  \times \left[ h (Y) \cup \cdots \cup h^m (Y) \right]
 \cup  \{- m (j - 1) \}
  \times \left[ Y \cup \cdots \cup h^{-m + 1} (Y) \right]  \\
  & \S \Z \times X.
\end{align*}
Then the $S_j$ are compact $\Z \times X$-sets, $s (S_j) = s (K)$, and
\begin{align*}
r (S_j) & =
  \left[ h^{m (j - 1) + 1} (Y) \cup h^{m (j - 1) + 2} (Y)
                     \cup \cdots \cup h^{m j} (Y) \right]  \\
   & \hspace*{3em} \mbox{}  \cup
  \left[ h^{- m (j - 1)} (Y) \cup h^{-[ m (j - 1) + 1]} (Y)
              \cup \cdots \cup h^{-(m j - 1)} (Y) \right]  \\
 & = X_1 \times \left[
     \{ h_2^{m (j - 1) + 1} (z), \, \dots, \, h_2^{m j} (z) \} \cup
        \{ h_2^{- m (j - 1)} (z), \, \dots, \, h_2^{- (m j - 1)} (z) \}
    \right].
\end{align*}
Since the \hme\  $h_2$ has no periodic points,
these sets are pairwise disjoint.
This completes the proof that if $K \subset G \setminus G_0$ is compact,
then $s (K)$ is thin in $G_0$,
and hence also the proof of condition~(1) of Definition~\ref{TAFGDfn}.

The reduced \ca\  $C_{\mathrm{r}}^* (\Z \times X)$ is isomorphic
to the \tgca\  $C_{\mathrm{r}}^* (\Z, X, h) = C^* (\Z, X, h)$,
by Proposition~\ref{CrPisGpC}.
We can use Theorem~3.1 of~\cite{Pt2} to show that $C^* (\Z, X, h)$
does not have stable rank one, but it is just as easy to do this
directly.
Let $u \in C^* (\Z, X, h)$ be the canonical unitary, satisfying
$u f u^* = f \circ h^{-1}$ for $f \in C (X)$.
Set
\[
R = \{ 0, 1, \dots, \infty \} \times X_2 \subset X
\andeqn
p = \ch_R \in C (X) \subset C^* (\Z, X, h).
\]
Then one checks that $s = u p + (1 - p)$ satisfies
\[
s^* s = 1 \andeqn s s^* = 1 - \ch_{ \{ 0 \} \times X_2 } \neq 1.
\]
So $s$ is a nonunitary isometry in $C^* (\Z, X, h)$.
It follows that $s$ is not a norm limit of invertible elements of
$C^* (\Z, X, h)$, so this algebra does not have stable rank one.
(Compare with Examples~4.13 of~\cite{Rf}.)
\end{exa}

\end{document}